\documentclass{article}
\usepackage[utf8]{inputenc}

\usepackage{bbm}
\usepackage[giveninits=true,sorting=none]{biblatex}
\usepackage{hyperref}
\usepackage{graphicx}
\usepackage{subcaption}
\usepackage{enumerate}
\usepackage{booktabs}
\usepackage{dirtytalk}
\usepackage{sidecap}

\usepackage{xcolor}
\usepackage{amsmath}
\usepackage{amsfonts}
\usepackage{stmaryrd}
\usepackage{array}
\usepackage{caption}
\newcolumntype{M}[1]{>{\centering\arraybackslash}m{#1}}

\usepackage[ruled,vlined]{algorithm2e}

\title{Machine Learning based refinement strategies\\ for polyhedral grids with applications to\\ Virtual Element and polyhedral Discontinuous Galerkin methods} 

\author{P. F. Antonietti\thanks{MOX, Department of Mathematics, Politecnico di Milano, p.zza Leonardo da Vinci 32, I-20133, Milano, Italy (paola.antonietti@polimi.it, enrico.manuzzi@polimi.it).}, F. Dassi\thanks{Department of Mathematics and Applications, University of Milano-Bicocca, Via Cozzi 53, I-20153, Milano, Italy (franco.dassi@unimib.it)}, E. Manuzzi\footnotemark[1]}

\begin{document}

\maketitle
\section*{Abstract}
We propose two new strategies based on Machine Learning techniques to handle polyhedral grid refinement, to be possibly employed within an adaptive framework. The first one employs the k-means clustering algorithm to partition the points of the polyhedron to be refined. This strategy is a variation of the well known Centroidal Voronoi Tessellation. The second one employs  Convolutional Neural Networks to classify the \textit{“shape”} of an element so that “ad-hoc” refinement criteria can be defined. This strategy can be used to enhance existing refinement strategies, including the k-means strategy, at a low online computational cost. We test the proposed algorithms considering two families of finite element methods that support arbitrarily shaped polyhedral elements, namely the Virtual Element Method (VEM) and the Polygonal Discontinuous Galerkin (PolyDG) method. We demonstrate that these strategies do preserve the structure and the quality of the underlaying grids, reducing the overall computational cost and mesh complexity.\\
\newline
\textbf{Keywords:} polyhedral grid refinement, Machine Learning, Convolutional Neural Networks, k-means, Polyhedral Discontinuous Galerkin, Virtual Element Method.\\
\newline
\textbf{Abbreviations:} Machine Learning (ML), Centrodial Voronoi Tesselation (CVT), Convolutional Neural Networks (CNNs), Finite Element Methods (FEMs),Virtual Element Method (VEM), Polyhedral Discontinuous Galerkin (PolyDG).

\section{Introduction}
{\color{black} Many applications in the fields of Engineering and Applied Sciences, such as fluid-structure interaction problems, flow in fractured porous media, crack and wave propagation problems, are characterized by a strong complexity of the physical domain, possibly involving moving geometries, heterogeneous media, immersed interfaces (such as e.g. fractures) and complex topographies. Whenever classical Finite Element Methods (FEMs) are employed to discretize the underlying differential model, the process of grid generation can be the bottleneck of the whole simulation, as computational meshes can be composed only of tetrahedral, hexahedral, or prismatic elements. To overcome this limitation, in the last years there has been a great interest in developing FEMs that can employ general polygons and polyhedra as grid elements for the numerical discretizations of partial differential equations.
We mention the mimetic finite difference method \cite{hyman1997numerical,brezzi2005family,brezzi2005convergence,da2014mimetic}, the hybridizable discontinuous Galerkin method \cite{cockburn2008superconvergent,cockburn2009superconvergent,cockburn2009unified,cockburn2010projection}, the Polyhedral Discontinuous Galerkin (PolyDG) method \cite{hesthaven2007nodal,bassi2012flexibility,antonietti2013hp,cangiani2014hp,antonietti2016review,cangiani2017hp,antonietti2021high}, the Virtual Element Method (VEM) \cite{beirao2013basic,beirao2014hitchhiker,beirao2016virtual,da2016mixed,beirao2021recent,book.vem.sema.simai.2022} and the Hybrid High-Order method \cite{di2014arbitrary,di2015hybrid,di2015hybrid2,di2016review,di2019hybrid}.
This calls for the need to develop effective algorithms to handle polygonal and polyhedral grids and to assess their quality (see e.g. \cite{attene2019benchmark}).
For a comprehensive overview we refer to the monographs and special issues \cite{da2014mimetic,cangiani2017hp,di2016review,di2021polyhedral,beirao2021recent,book.vem.sema.simai.2022}.\\
Among the open problems, there is the issue of efficiently handling polytopic mesh refinement \cite{lai2016recursive,hoshina2018simple,berrone2021refinement}, i.e., partitioning mesh elements into smaller elements to produce a finer grid, and agglomeration strategies, i.e., merging mesh elements to obtain coarser grids \cite{chan1998agglomeration,antonietti2020agglomeration,bassi2012flexibility}. Indeed, as for standard triangular and quadrilateral meshes, during either refinement or agglomeration it is important to preserve the quality of the underlying mesh, since this might affect the overall performance of the method in terms of stability and accuracy. Indeed a suitable adapted mesh may allow to achieve the same accuracy with a much smaller number of degrees of freedom when solving the numerical problem, hence saving memory and computational power. However, since in such a general framework mesh elements may have any shape, there are not well established strategies to achieve effective refinement or agglomeration with a fast, robust and simple approach, contrary to classical tetrahedral, hexahedral and prismatic meshes. Moreover, grid agglomeration is a topic quite unexplored, because it is not possible to develop such kind of strategies within the framework of classical FEMs.\\}
\newline
In recent years there has been a great development of Machine Learning (ML) algorithms to enhance and accelerate numerical methods for scientific computing. Examples include, but are not limited to, \cite{raissi2019physics,raissi2018hidden,regazzoni2019machine,regazzoni2020machine,hesthaven2018non,ray2018artificial,antonietti2021accelerating,regazzoni2021machine}. In this work, we propose different strategies to handle polyhedral grid refinement, that exploit ML techniques to extract information about the “shape” of mesh elements in order to refine them accordingly. Starting from the approach developed in \cite{ANTONIETTI2022110900} in two dimensions, here we address the three-dimensional case. In particular, we employ the k-means clustering algorithm \cite{hartigan1979algorithm,likas2003global} to partition the points of the polyhedron to be refined. Learning a clustered representation allows to better handle situations where elements have no particular structure, while preserving the quality of the grid. This strategy is a variation of the well known Centroidal Voronoi Tessellation \cite{du1999centroidal,liu2009centroidal}. We also explore how to enhance existing refinement strategies, including the k-means, using Convolutional Neural Networks (CNNs), powerful function approximators successfully employed in many areas of ML, especially computer vision \cite{lecun2015deep,albawi2017understanding}. We show that CNNs can be employed to identify correctly the “shape” of a polyhedron, without resorting to the explicit calculation of geometric properties which may be expensive in three dimensions. This information can then be exploited to design tailored strategies for different families of polyhedra. According to the classification, the most suitable refinement strategy is then selected among the available ones, therefore enhancing performance. The proposed approach has several advantages:
\begin{itemize}
    \item it helps preserving structure and quality of the mesh, since it can be easily tailored for different types of elements;
    \item it can be combined with suitable user-defined refinement criteria, including the k-means strategy;
    \item it is independent of the differential model at hand and of the numerical method employed for the discretization;
    \item the overall computational cost is kept low, since CNNs can be trained offline once and for all.
\end{itemize}
To investigate the capabilities of the proposed approaches, we consider a second-order model problem discretized by either the VEM and the PolyDG method. We measure effectiveness through an analysis of quality metrics and accuracy of the discretization error.\\
\newline
The paper is organized as follows. In Section 2 we present the k-means refinement strategy and other possible refinement criteria for polyhedra. In Section 3 we propose a general framework to enhance existing refinement strategies using CNNs. In Section 4 we validate the proposed refinement strategies on a set of polyhedral meshes and we measure their effectiveness using the quality metrics introduced in \cite{attene2019benchmark}. In Section~5 we
present some computations obtained with a uniform and a priori adaptive mesh refinement process for a differential model approximated either by VEMs and PolyDG methods. In Section 6 we draw some conclusions.

\section{Refinement Strategies}
\begin{figure}
    \centering
    \includegraphics[width = 0.8\linewidth]{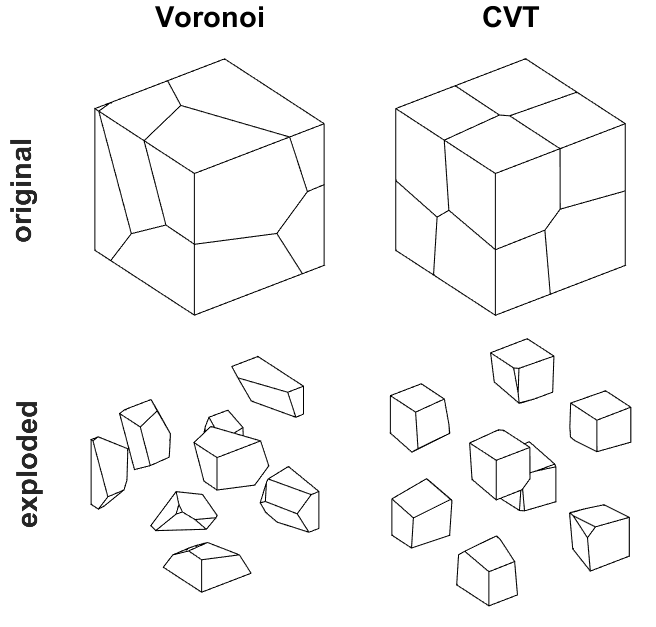}
    \caption{A Voronoi tessellation (top left) and a CVT (top right) of a cube and their exploded versions (bottom). “Small" edges and faces are generated.}
    \label{fig:voronoi}
\end{figure}

A polyhedral mesh can be generated, e.g., by standard triangulation algorithms or by hexahedral elements in the structured case. Another possibility is to use a Voronoi Tessellation \cite{hoshina2018simple}: suitable points, called seeds, are chosen inside the polyhedron and each element of the new partition is the set of points which are closer to a specific seed (see Figure~\ref{fig:voronoi} left column).
A Centroidal Voronoi Tessellation (CVT) \cite{du1999centroidal,liu2009centroidal} is a Voronoi tessellation whose seeds are the centroids (centers of mass) of the corresponding new elements (see Figure~\ref{fig:voronoi} right column). This strategy produces elements which are shape-regular and of similar size (these properties will be more formally addressed in Section \ref{section quality}). However, it has a considerable computational cost and it may produce elements with many small edges and faces, which implies higher costs in terms of memory storage and algorithm complexity. Moreover, although small edges and faces do not necessarily deteriorate the accuracy of numerical methods, they are in general not beneficial \cite{brenner2018virtual,beirao2017stability,droniou2021robust,mu2015shape}. As we will see in the following, these algorithms for grids generation can be adapted to perform mesh refinement.\\
\newline
In order to refine a three dimensional polyhedron a possible strategy is to subdivide the polyhedron along a chosen preferential direction \cite{berrone2021refinement}. It has a low computational cost and allows to adjust the cutting plane at each iteration. For instance it is possible to choose such plane so that we avoid small edges and consequently yield simpler mesh elements. In particular, we propose the general refinement strategy described in Algorithm \ref{alg:general refinement strategy}, where the \say{size} or \say{diameter} of a polyhedron $P \subset \mathbb{R}^3$ is defined, as usual, as
$$\textrm{diam}(P) = \sup \{||x-y||,\ x,y \in P\},$$
where $||\cdot||$ is the standard Euclidean distance. We also recall that given a polyhedral mesh, i.e., a set of non-overlapping polyhedra $\{P_i\}_{i=1}^{N_P}$, $N_P \geq 1$ that covers a domain $\Omega$, the mesh size is defined as
$$h = \max_{i=1,..., N_P} \textrm{diam}(P_i).$$
Possible ways of choosing the cutting plane in step 1 will be discussed in Section \ref{section cut plane}.
\begin{algorithm}
\SetAlgoLined
\LinesNumbered
\KwIn{polyhedral element $P$ to refine, plane tolerance \texttt{tol}, maximum number of elements \texttt{nmax}, target size $\Bar{h}$.}
\KwOut{partition of $P$ into polyhedral sub-elements.}
\DontPrintSemicolon
\BlankLine
\BlankLine
Choose a plane along which to slice the element. \;
If there are vertices which closer than \texttt{tol} to the plane, adjust the plane so that it passes from them. \;
If the plane passes the \say{validity check}, slice the element into two sub-elements, otherwise use the \say{emergency strategy}.\;
Until there are less than \texttt{nmax} elements, repeat the above steps for each element with a size above $\Bar{h}$. \;
\caption{General refinement strategy of a polyhedron $P$}
\label{alg:general refinement strategy}
\end{algorithm}
The \say{validity check} in step 3, to see if it is possible to slice the element along the chosen plane, may encompass several requirements.

\begin{figure}
    \centering
    \includegraphics[width = 0.75\linewidth]{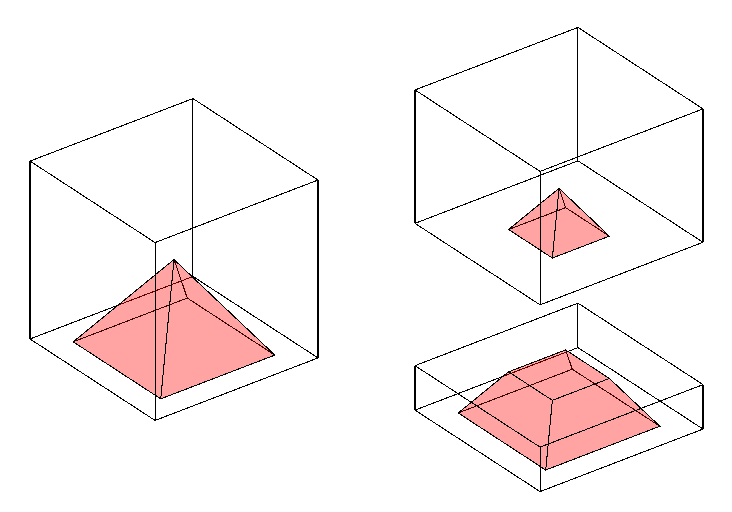}
    \caption{Left: a non-convex polyhedron is obtained from a cube (transparent) by removing a pyramid (red). Right: the polyhedron is sliced along a plane generating two polyhedra, one of which has a hole.}
    \label{fig:degenerate}
\end{figure}

    \paragraph{Geometrical.} Avoid generating elements with holes and other degenerations, which may occur if the polyhedron is non-convex (see Figure \ref{fig:degenerate}).
    \paragraph{Numerical.} Avoid generating small edges and faces according to a prescribed threshold, which may still occur because in step 2 a plane can be adjusted to pass from 3 points.\\[0.5em]
If the chosen plane does not pass the \say{validity check} a possible \say{emergency strategy} is the following: consider small perturbations of position and orientation of the original plane, until a valid configuration is found.

\subsection{Choosing the cutting direction} \label{section cut plane}
\begin{figure}
    \centering
    \includegraphics[width = \linewidth]{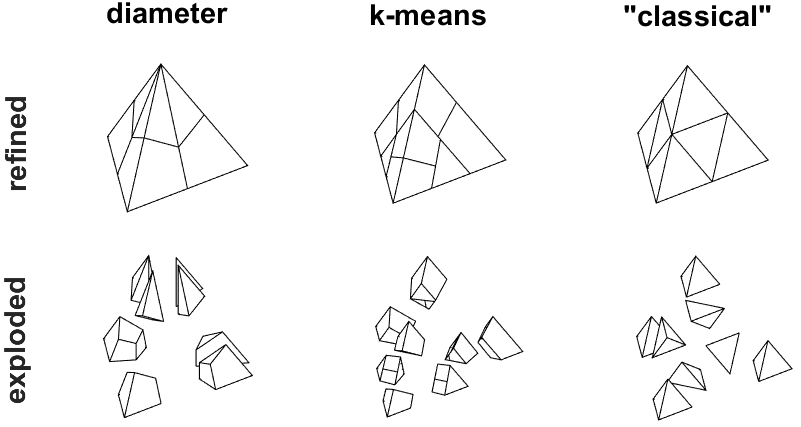}
    \caption{A tetrahedron is refined using diameter, k-means and {\color{black}``classical"} strategies (top row) and their exploded versions (bottom row).}
    \label{fig:tet comparison}
\end{figure}

Several strategies can be designed by making different choices of the cutting direction in step 1 of Algorithm \ref{alg:general refinement strategy}. In Figure \ref{fig:tet comparison}  the proposed strategies are applied on a tetrahedron.
Notice that which strategy is the most effective may depend on several factors, such as the initial geometry of the mesh, the quality of the refined elements, the computational cost, the stability of the numerical scheme, or the regularity of the solution.
{\color{black}
    \paragraph{Diameter strategy.} The diameter of a polyhedron is identified by the two farthest apart vertices of the element, say $v_1$ and $v_2$. The cutting plane has origin $x_0 = (v_1+v_2)/2$ and is orthogonal to the direction of the diameter, i.e., it has normal $n = (v_2-v_1)/||v_2-v_1||$. This is a greedy algorithm to halve the size of the element. It has the advantage of easily computing the cutting plane, but it may produce skewed elements because the overall shape of the polyhedron is not taken into account.

    \paragraph{K-means strategy.} Points belonging to the polyhedron are grouped into two clusters using the k-means algorithm \cite{hartigan1979algorithm,likas2003global} and which are then separated by the cutting plane, as described in Algorithm \ref{alg:k-means}. This is equivalent to a CVT with only two seeds. This strategy produces rounded elements of similar size, but it has a higher computational cost. In principle, different clustering algorithms may be used to learn different representations of the element, such as hierarchical clustering \cite{johnson1967hierarchical,murtagh2012algorithms} or self-organizing maps \cite{ritter1992neural,kangas1990variants}.
    }

\begin{algorithm}
{\color{black}
\SetAlgoLined
\LinesNumbered
\KwIn{polyhedron $P$ to refine, number of grid points $n$.}
\KwOut{cutting plane of $P$.}
\DontPrintSemicolon
\BlankLine
\BlankLine
Construct a grid of $n$ of points inside the smallest domain of the form $[x_{1},x_{2}] \times [y_{1},y_{2}] \times [z_{1},z_{2}]$ which contains $P$.\;
Remove the points of the grid outside $P$ \cite{lane1984efficient}, which can be done efficiently when $P$ is convex.\;
Randomly sample two points and label them centroids $c_1$ and $c_2$.\;
Assign each point to the cluster with the closest centroid in $L^2$ norm.\;
Compute the average of points in each cluster to obtain two new centroid $c_1$ and $c_2$.\;
Repeat steps 4 and 5 until cluster assignments do not change, or the maximum number of iterations is reached.\;
The cutting plane has normal $n = (c_2-c_1)/||c_2-c_1||$ and origin $x_0 = (c_1+c_2)/2$.\;
\caption{K-means cutting plane}
\label{alg:k-means}
}
\end{algorithm}

\begin{figure}
    \centering
    \includegraphics[width = 1\linewidth]{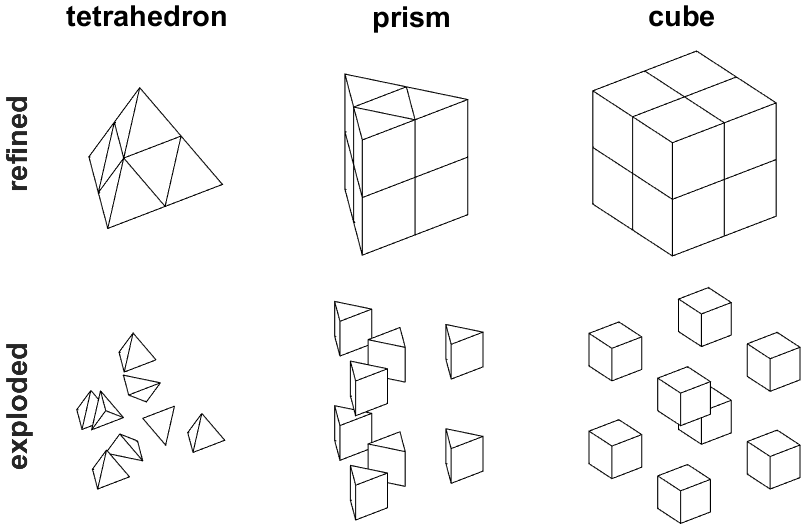}
    \caption{{\color{black}``Classical"} refinement strategies for a tetrahdron, a prism and a cube (top row) and their exploded version (bottom row), employed in classical FEMs.}
    \label{fig:fem shapes}
\end{figure}

    \paragraph{{\color{black}``Classical"} strategies.} If the shape of the polyhedron is known (up to rotation, scaling, stretching or reflection) it is possible to design tailored refinement strategies using multiple cutting planes. This is the case of tetrahedra, cubes and prisms, which are employed in classical FEMs. These refinement strategies, shown in Figure \ref{fig:fem shapes}, produce regular elements at low computational cost. Moreover, since the new elements have the same shape of the original one, the algorithm can be applied recursively.
    Strategies for other polyhedral element can be designed. The main drawback is that these strategies cannot be applied if the shape of the polyhedron is unknown.

\section{Enhancing refinement strategies using CNNs}
Consider elements at the top row of Figure \ref{fig:why CNN} from left to right:

\begin{figure}
    \centering
    \includegraphics[width = \linewidth]{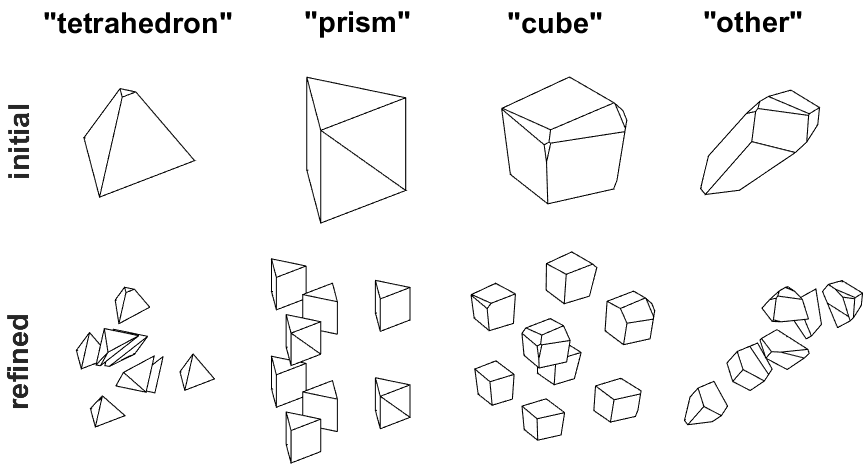}
    \caption{Top row: four polyhedra are classified as ``tetrahedron", ``prism", ``cube" and ``other", according to their overall shape and not to their geometrical description. Bottom row: the polyhedra are refined with the corresponding ``classical" strategies according to their labels (columns 1-3) and the k-means strategy (column 4).}
    \label{fig:why CNN}
\end{figure}

\begin{enumerate}
    \item A tetrahedron where a corner has been cut. This situation may arise when slicing a background mesh, e.g., in order to model a fracture or a soil discontinuity.
    \item A prism where one of its faces is divided into two, due to the presence of non-matching grids. This situation can arise each time a neighbouring mesh element is refined, because this may require to partition a face which is shared by two elements.
    \item An element of a CVT. These tessellations are employed for meshing irregular domains with few regular elements.
    \item A skewed element with no particular structure. These elements can be found in standard Voronoi tessellations, or may appear even after refining a regular element.
\end{enumerate}
Despite the fact that the geometrical representations of elements 1-3 in Figure \ref{fig:why CNN} do not match any classical polyhedron, their overall shapes are similar to a tetrahedron, a prism and a cube, respectively. Therefore refining them with the corresponding shape strategies would produce regular elements at low computational cost. This can be done by computing cutting directions on a \say{reference shape} which can be constructed using Algorithm \ref{alg:reference shape}.
\begin{algorithm}
\SetAlgoLined
\LinesNumbered
\KwIn{polyhedron $P$, shape type $S$ (tetrahedron, cube, prism, ...).}
\KwOut{reference shape of $P$.}
\DontPrintSemicolon
\BlankLine
\BlankLine
Find $n$ vertices of the polyhedron far apart from each others and from the centroid of the element, using the farthest first traversal algorithm \cite{le2020farthest}, where $n$ is the number of vertices of $S$.\;
Apply the convex hull algorithm \cite{barber1996quickhull} to the $n$ vertices to obtain a triangulated surface.\;
Merge the $m$ couples of neighbouring triangles of the surface whose angle is closest to $180^\circ$, where $m$ is the number of quadrilateral faces of $S$.\;
\caption{Computing the reference shape}
\label{alg:reference shape}
\end{algorithm}
The last element does not have a particular shape, and therefore we apply the k-means strategy, which is more robust and flexible even though more expensive. The refined elements are shown at the bottom row of Figure \ref{fig:why CNN}. However, in order to apply the correct strategy for each element, we need an algorithm to access whether the shape of a polyhedron is more similar to a tetrahedron, a prism, a cube or none of these. We will use labels \say{tetrahedron}, \say{prism}, \say{cube} and \say{other}, respectively, and we will refer to such labels as \say{equivalence classes}. More in general, given a set of polyhedra to refine, $\mathcal{P} = \{P_1,P_2, ...\}$, and a set of possible refinement strategies, $\mathcal{R} = \{R_1,R_2,...,R_n\}$, we want to find a classifier $F:\mathcal{P} \rightarrow \mathcal{R}$ that assigns to every polyhedron the most effective refinement strategy in terms of quality of the refined elements and computational cost. The effectiveness will be measured with suitable metrics later in Section \ref{section quality}. 

\subsection{Polyhedra classification using CNNs} 
\begin{figure}
    \centering
    \includegraphics[width = 0.45\linewidth]{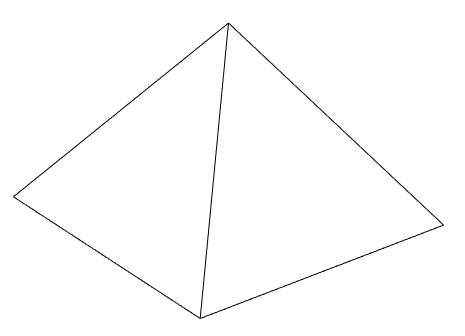}
    \hfill
    \includegraphics[width = 0.45\linewidth]{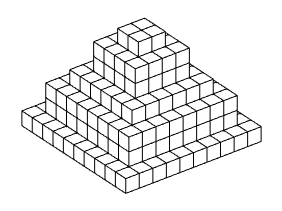}
    \caption{A pyramid (left) and its three-dimensional image (right). Pixels have value 1 inside the polyhedron (represented with cubes) and 0 outside (not represented).}
    \label{fig:3D image}
\end{figure}

In principle, any classifier may be used for polyhedra. However, considering geometrical properties of the element alone, such as connectivity and area of the faces, may be too complex to operate a suitable classification: for example, the element labeled as \say{tetrahedron} in Figure \ref{fig:why CNN} is actually a deformed prism. Approaches like Shape Analysis \cite{dryden2016statistical}, which relies on applying transformations the polyhedron in order to match a reference shape, are too expensive in three-dimensions and seem to lack of the required flexibility: for example there is no reference shape for class \say{other}. Instead, providing samples of polyhedra together with the desired label can be easily done: for classes ``tetrahedron", ``prism" and ``cube"  we start from the corresponding regular polyhedra and then apply small deformations and rotations, while for class ``other" we consider polyhedra from several Voronoi tessellations. This database can be used to “train” our function, i.e., to tune its parameters in order to obtain the desired classifier. This framework, which exploits labeled data, is known as “supervised learning” and CNNs are powerful function approximators in the context of image classification. Indeed, the most intuitive way to asses the shape of an object is by visual inspection and polyhedra can be easily converted into binary images: each pixel assumes value 1 inside the polyhedron and 0 outside \cite{lane1984efficient} (see Figure \ref{fig:3D image}), which can be done efficiently in the convex case.
Mathematically, while the geometrical representation of a polyhedron is given by its boundary, the information about the overall shape is given by the distribution of its volume, i.e., the location of its pixels. Notice that a three-dimensional image contains the same information (up to a scaling and a translation) required by the k-means algorithm, which is a form of \say{unsupervised} learning because only unlabeled data are provided. The proposed CNN-enhanced refinement strategy is describe in Algorithm \ref{alg:CNN ref}.
\begin{algorithm}
\SetAlgoLined
\LinesNumbered
\KwIn{polyhedron $P$ to refine.}
\KwOut{partition of $P$ into polyhedral sub-elements.}
\DontPrintSemicolon
\BlankLine
\BlankLine
Construct a binary image of $P$.\;
Classify the image using a CNN and obtain the label $\mathcal{L}$.\;
\BlankLine
\If{$\mathcal{L} =$ ``tetrahedron", ``prism" or ``cube"}{
Apply Algorithm \ref{alg:reference shape} to compute the corresponding reference shape.\;
Compute the cutting directions on the reference shape using the corresponding ``classical" strategy.\;
Refine $P$ by applying the cutting directions using Algorithm \ref{alg:general refinement strategy}.\;
}
\BlankLine
\If{$\mathcal{L} =$ ``other"}{
Refine $P$ using the k-means strategy, i.e., computing the cutting directions using Algorithm \ref{alg:k-means} and then using Algorithm \ref{alg:general refinement strategy}.\;
}

\caption{CNN-enhanced refinement strategy}
\label{alg:CNN ref}
\end{algorithm}
In general, the use of CNNs should be seen as tool to enhance existing refinement strategies, not as a refinement strategy on its own in competition with the others. In the case of ``classical" strategies, the CNN is extending their domain of applicability to general  polyhedra, which would not be possible otherwise.

{\color{black}
\subsection{Supervised learning for image classification}
Consider a three dimensional binary image $\textbf{B} \in \{0,1\}^{n \times n \times n}$, $n\geq 1,$ and the corresponding label vector $y = ([y]_j)_{j =1,..,\ell} \in  [0,1]^\ell,$
where $\ell \geq 2$ is the total number of classes, and $[y]_j$ is the probability of $\textbf{B}$ to belong to the class $j$ for $j = 1:\ell$. In a supervised learning framework, we are given a dataset of desired input-output couples $\{(\textbf{B}_i,y_i)\}_{i=1}^N$, where $N$ is number of labelled data.
We consider then an image classifier represented by a function of the form $F: \mathbb{R}^{n \times n \times n} \rightarrow~(0,1)^\ell,$ in our case a CNN, parameterized by $w \in \mathbb{R}^M$ where $M\geq 1$ is the number of parameters. Our goal is to tune $w$ so that $F$ minimizes the data misfit, i.e.,
$$\min_{w \in \mathbb{R}^M} \sum_{i \in I} l(F(\textbf{B}_i),y_i),$$   
where $I$ is a subset of $\{1,2,...,N\}$ and $l$ is the cross-entropy loss function defined~as
$$ l(F(\textbf{B}),y) = \sum_{j = 1}^\ell -[y]_j\log[F(\textbf{B})]_j.$$
This optimization phase is also called ``learning" or ``training phase". During this phase, a known shortcoming is ``overfitting": the model fits very well the data used in the training phase, but performs poorly on new data. We proceed in a standard way \cite{bishop2006pattern}, i.e., we split the data into
\begin{itemize}
    \item \emph{training set}: used to tune the parameters during the training phase;
    \item \emph{validation set}: used to monitor the model capabilities during the training phase, on data that are not in the training set. The training is halted if the error on the validation set starts to increase;
    \item \emph{test set}: used to access the actual model performance on new data after the training.
\end{itemize}
While the training phase can be computationally demanding, because of the large amount of data and parameters to tune, it needs to be performed offline once and for all. Instead, classifying online a new image using a pre-trained model is computationally fast: it requires only to evaluate $F$ on a new input. The predicted label is the one with the highest estimated probability.}
{\color{black}
\subsection{CNN architecture}
CNNs are parameterized functions, in our case of the form $$\text{CNN}:~\mathbb{R}^{n \times n \times n}~\rightarrow (0,1)^\ell,\hspace{0.5cm} n,\ell \geq 2,$$ constructed by composition of simpler functions called ``layers of neurons" \cite{lecun2015deep}. We basically use these types of layers.

    \paragraph{Convolutional layers.} Linear maps performing a convolution of the input, which can be viewed as a filter scanning through the image, extracting local features that depend only on small subregions of the image. This is effective because a key property of images is that close pixels are more strongly correlated than distant ones. The scanning filter mechanism provides the basis for the invariance of the output to translations and distortions of the input image \cite{bishop2006pattern}. Multiple convolutions may be performed in a single layer, each of which is called a \say{feature map}.

    \paragraph{Pooling layers.} Maps used to perform down-sampling by merging semantically similar features. They perform operations region-wise, such as averaging or computing the maximum. These regions are rectangular and they are computed by fixing a size and then moving along the image of an amount of pixels called \say{stride}. They improve the invariance of the output with respect to translations of the input.
    
   \paragraph{Activation functions.} Maps used to introduce non-linearity, such as the rectified linear unit $$[\textsc{ReLU}(x)]_i = \max(0,x_i).$$

    \paragraph{Dense layers.} Generic linear maps of the form $\textsc{Linear}:  \mathbb{R}^{m} \rightarrow \mathbb{R}^\ell,$ $m,\ell \geq 1$ defined by parameters to be tuned. They are used to separate image features extracted in the previous layers.
    
    \paragraph{Softmax.} We define the function $\textsc{Softmax}: \mathbb{R}^\ell \rightarrow (0,1)^\ell,$ where $\ell\geq 2$ is the number output classes, $$[\textsc{Softmax}(x)]_i = \frac{e^{x_i}}{\sum_{j=1}^\ell e^{x_j}}.$$ They are used to assign a probability to each class.\\[0.5em]    
In practise, subsequent application of convolutional, activation and pooling layers may be used to obtain a larger degree of invariance to input transformations such as rotations, distortions, etc...}


\subsection{CNN training}
\begin{figure}
    \centering
    \includegraphics[width = 0.75\linewidth]{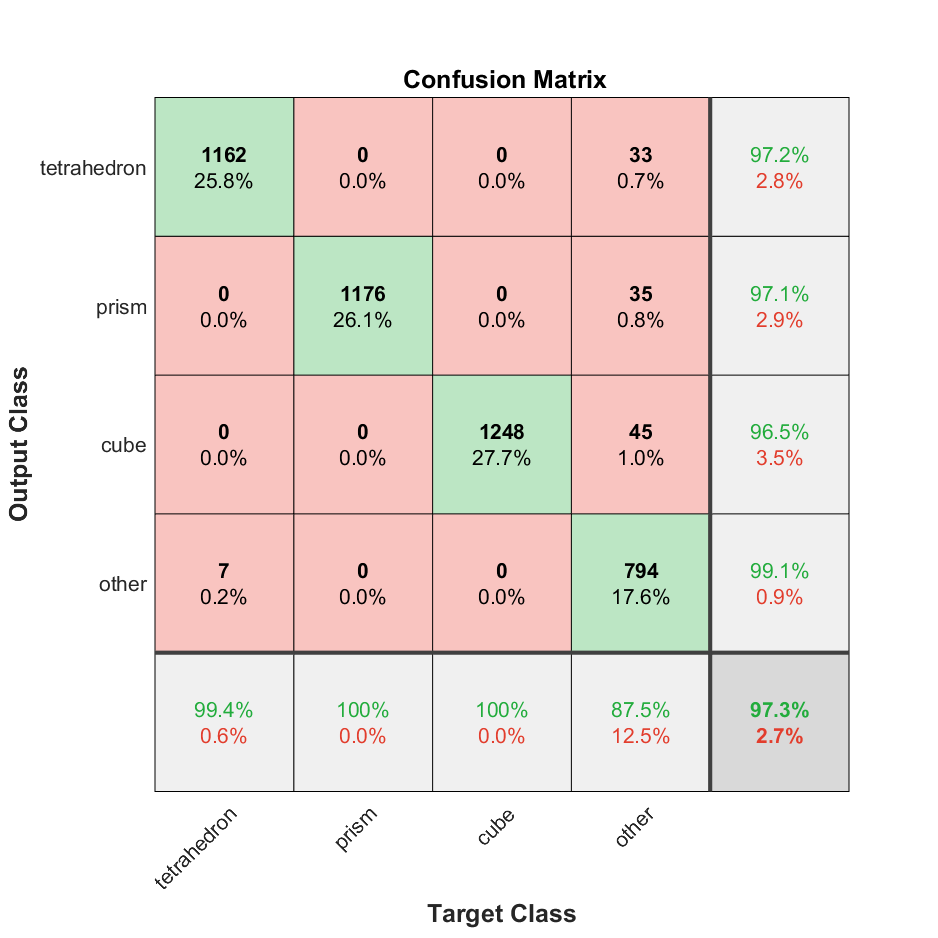}
   \caption{Confusion matrix for polyhedra classification using a CNN. Class “other" is reasonably the most misclassified, because it includes polyhedra with very different shapes, some of which are actually deformed tetrahedra, prisms and cubes.}
    \label{fig:confusion matrix}
\end{figure}

We have used a database of 22500 binary images of size $16 \times 16 \times 16$ pixels, with labels ``tetrahedron", ``prism", ``cube" and ``other". The data were divided into  training, validation and test set with ratios 60\%-20\%-20\% respectively. They were generated as follows:
\begin{itemize}
    \item for classes ``tetrahedron", ``prism" and ``cube"  we generated 6000 images each, starting from the corresponding regular polyhedra, and then randomly applying perturbations, such as rotation, scaling, stretching and reflection;
    \item for the class ``other" we generated 4500 images from elements of several Voronoi tessellations, where seeds have been randomly located.
\end{itemize}
We generated less images for class ``other", because among all of the randomly generated Voronoi elements, some of them will be actually deformed tetrahedra, prisms and cubes. In these doubtful situations, the CNN will be more likely to classify polyhedra as having a particular shape rather assigning them label ``other". Indeed, the CNN is capable to learn the whole data distribution, including the relative frequency of each class. In this way, we resort to the k-means only when really needed.\\
Generating data which are representative of the application of interest is the most critical part of the process: the CNN can easily misclassify new samples if they are too different from the training set, even if the accuracy on the test set is extremely high. In three dimensions, the variability of polyhedra is high and generating representative sample is much more complex than in the two-dimensional case \cite{ANTONIETTI2022110900}. It is important to have a class ``other" to which unknown polyhedra will be assigned, in order to treat them with a robust backup strategy, such as the k-means. This should not discourage from using ``classical" strategies, because when employed they are likely to produce regular elements.\\
\newline
The CNN architecture we used is $\text{CNN}: \{0,1\}^{16 \times 16 \times 16} \rightarrow (0,1)^4,$
where
\begin{align*} 
\text{CNN} =\ 
& \textsc{Conv}(\Bar{f}  = 8,m =8)   \rightarrow \textsc{ReLU}\rightarrow \textsc{Pool}(\Bar{f}  = 2, s = 2) \rightarrow \\
&  \textsc{Conv}(\Bar{f}  = 4, m =8)  \rightarrow \textsc{ReLU} \rightarrow  \textsc{Pool}(\Bar{f}  = 2, s = 2)\rightarrow \\ 
& \textsc{Conv}(\Bar{f}  = 2,m =8)   \rightarrow \textsc{ReLU}\rightarrow \textsc{Pool}(\Bar{f}  = 2, s = 2) \rightarrow \\
& \textsc{Linear} \rightarrow \textsc{Softmax},
\end{align*}
where $\textsc{Conv}(\Bar{f} , s)$ is convolutional layer with filter size $\Bar{f} $ (i.e. a window $\Bar{f} \times \Bar{f}  \times \Bar{f} $) and $m$ feature maps, $\textsc{Pool}(\Bar{f} , s)$ is an avergage pooling layer with filter size $\Bar{f} $ and stride $s$,  $\textsc{ReLU}$ and $\textsc{Linear}$ were defined the previous Section 2.2, and for any two compatible functions $f_1$ and $f_2$ the arrow notation is defined as $f_1 \rightarrow f_2 = f_2 \circ f_1$. We employed average pooling layers, rather than max pooling layers, because they have the effect to blurry and coarse-grain image features. This has been done in order to focus the CNN on the overall shape of the polyhedron, and not on small scale details and corners. Three sequences of $\textsc{Conv}-\textsc{ReLU}-\textsc{Pool}$ layers seem to be sufficient to enforce invariance with respect to rotations of the input polyhedron. The size of the images is  $16\times 16 \times 16$ pixels, i.e., a resolution large enough to apply 3 pooling layers, whose effect is to halve the size of the image across all the dimensions. Such a small resolution allows to quickly generate and classify the image of a polyhedron, which is important in an online setting. Moreover, a smaller CNN is more likely to be robust with respect to samples very different from the training data. It was empirically observed that a larger resolution was not needed for the considered dataset. This is motivated by the fact that the approximation properties of neural networks generally depend on the dimension of the data manifold, not on the dimension of the input space \cite{petersen2020neural}. A larger resolution may be needed to classify polyhedra characterized by smaller scale features, e.g., to correctly distinguish between a dodecahedron and an icosahedron. As shown in Figure \ref{fig:confusion matrix}, the class “other" is reasonably the most misclassified, because it includes polyhedra with very different shapes, some of which are actually deformed tetrahedra, prisms and cubes.
Training the CNN using the Adam optimizer \cite{kingma2014adam} takes approximately 15 minutes with MATLAB2020b on a Windows OS 10 Pro 64-bit, CPU NVidia GeForce GTX 1050 Ti and 8244Mb Video RAM.

\section{Validation on a set of polyhedral grids} \label{section quality}
\begin{figure}[p!]
    \centering
    \includegraphics[width = \linewidth]{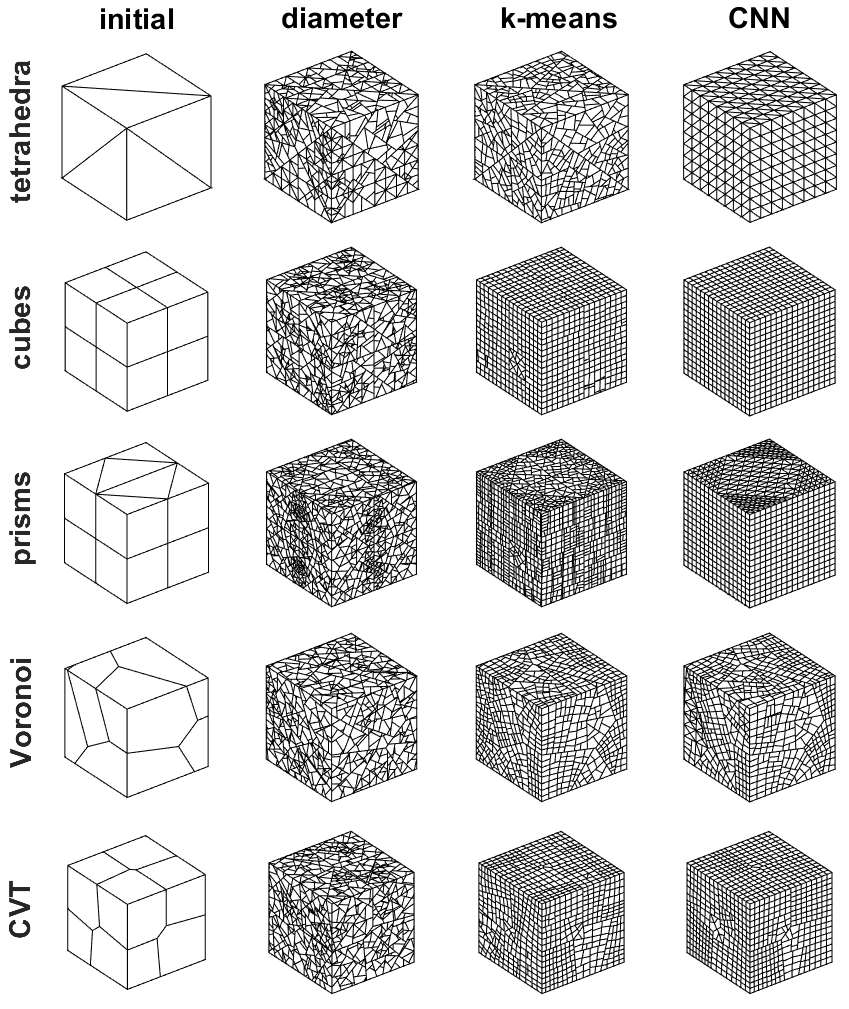}
  \caption{In the first column, coarse grids of domain $(0,1)^3$: a grid of tetrahedra, a grid of cubes, a grid of prisms, a Voronoi grid with random seeds location, and a CVT. Second to fourth columns: refined grids obtained after three steps of uniform refinement based on employing the k-means strategy (k-means) and the CNN-enhanced strategy (CNN). Each row corresponds to the same initial grid, while each column corresponds to the same refinement strategy.}
    \label{fig:corse grids refined}
\end{figure}

\begin{figure}
    \centering
    \includegraphics[width = \linewidth]{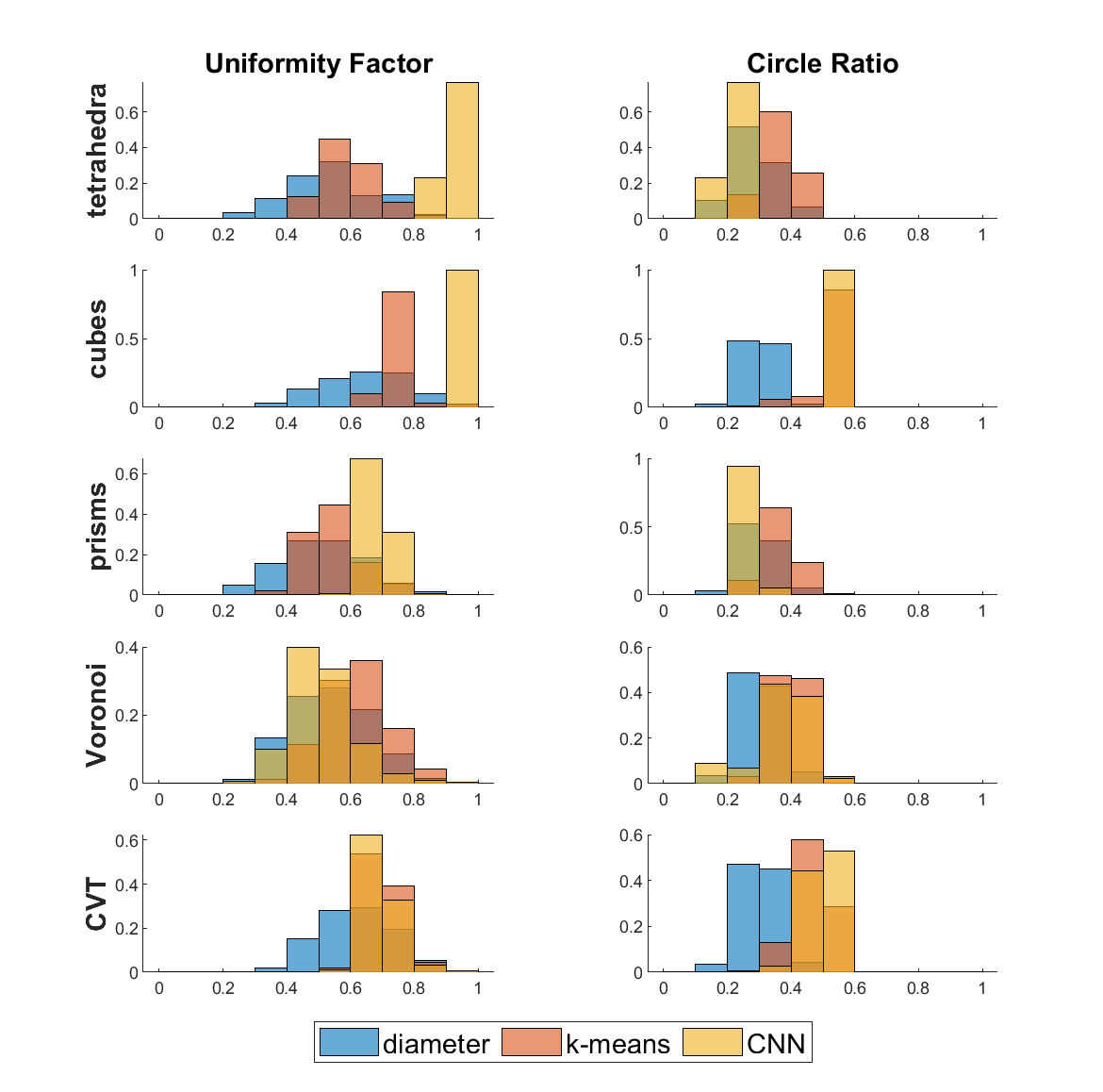}
  \caption{Histograms of the computed quality metrics (Uniformity Factor and Circle Ratio on the x-axis and percentage of grid elements on the y-axis) for the refined grids reported in Figure~\ref{fig:corse grids refined} (second to fourth column), obtained based on employing different refinement strategies (diameter, k-means, CNN).}
    \label{fig:quality metrics}
\end{figure}

\begin{figure}
    \centering
    \includegraphics[width = \linewidth]{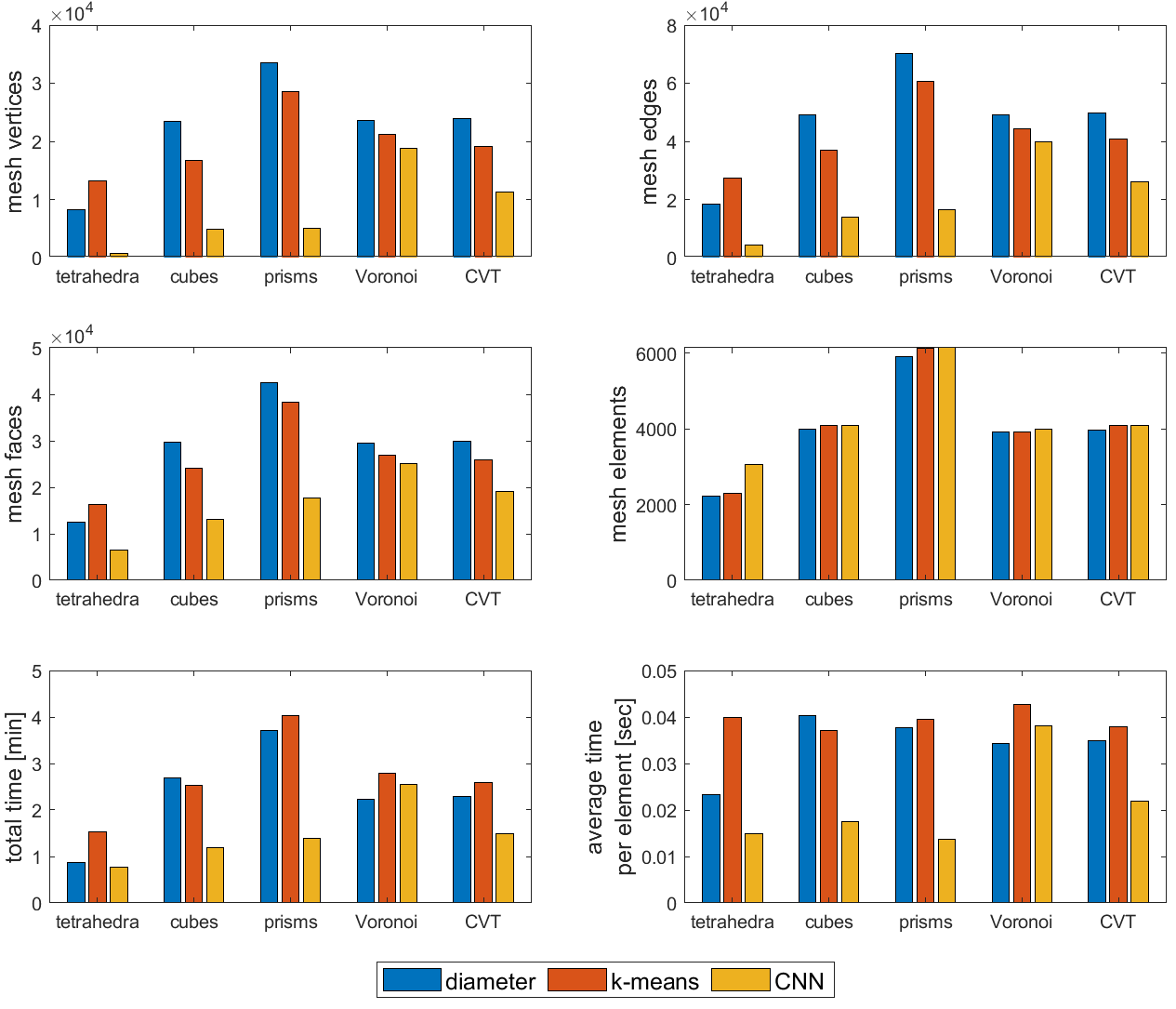}
  \caption{Statistics of computational complexity (number of vertices, edges, faces, elements, total refinement time and average refinement time per element) for the refined grids reported in Figure~\ref{fig:corse grids refined} (second to fourth column), obtained based on employing different refinement strategies (diameter, k-means, CNN).}
    \label{fig:complexity}
\end{figure}

In this section we compare the performance of the proposed algorithms. To evaluate the quality of the refined grids, we employ the following quality metrics introduced in \cite{attene2019benchmark}:
\begin{itemize}
    \item \textit{Uniformity Factor} (UF): ratio between the diameter of an element $P$ and the mesh size
    $$\textrm{UF}(P) = \frac{\textrm{diam}(P)}{h}.$$
    This metric takes values in $[0,1]$. The higher its value is the more mesh elements have comparable sizes.
    \item \textit{Circle Ratio} (CR): ratio between the radius of the inscribed circle and the radius of the circumscribed circle of a polyhedron $P$
    $$\textrm{CR}(P) = \frac{\max_{\{B(r) \subset P\}}  r }{\min_{\{P \subset B(r)\}} r },$$
	where $B(r)$ is a ball of radius $r$. For the practical purpose of measuring the roundness of an element the radius of the circumscribed circle has been approximated with $\textrm{diam}(P) / 2$. This metric takes values in $[0,1]$. The higher its value is the more rounded mesh elements are.
\end{itemize}
Moreover, to measure the computational complexity of the different refinement strategies, for the refined grids we will consider the number of vertices, of edges, of faces, of elements, the total refinement time and the average refinement time per element.\\
{\color{black} In the following, we apply Algorithm \ref{alg:general refinement strategy} to a polyhedron $P$ using a tolerance \texttt{tol}~$= \textrm{diam}(P) \cdot 10^{-3}$, a maximum number of elements \texttt{nmax} = 8 (needed to apply the \say{classical} refinement strategy for cubes), and a target size $\Bar{h} = \textrm{diam}(P)/2$. Concerning the k-means strategy, we apply Algorithm \ref{alg:k-means} using a grid of $n = 20^3$ points, which is comparable to the size of $16^3$ pixels required by the CNN. This value of $n$ may seem high, but it is needed in order compute the cutting plane with a satisfactory accuracy. Indeed, many grid points lie outside the polyhedron and therefore they are not included in the computation of the clusters, which is also the most expensive part of the process.} We consider five different coarse grids of domain $(0,1)^3$: a grid of tetrahedra, a grid of cubes, a grid of prisms, a Voronoi grid with random seeds location, and a CVT. In Figure~\ref{fig:corse grids refined} these grids have been refined uniformly, i.e., each mesh element has been refined, for three times using the diameter strategy (diameter), the k-means strategy (k-means) and the CNN-enhanced strategy (CNN).
As we can see, the initial mesh geometry seems to be disrupted by the diameter strategy, approximately preserved by the k-means strategy, and well preserved by the CNN-enhanced strategy.\\
In Figure \ref{fig:quality metrics} we show the computed quality metrics for the uniformly refined grids.
As we can see, in general quality metrics are lower for the diameter strategy. The k-means strategy produces elements with higher Circle Ratio, while the CNN-enhanced strategy with higher of Uniformity Factor, although in many cases both of them are comparable.\\
In Figure \ref{fig:complexity} we show, for all the considered grids, the number of vertices, edges, faces, elements, the total refinement time and the average refinement time per element. These results have been obtained by using MATLAB2020b on a Windows OS 10 Pro 64-bit, Intel(R) Core(TM) i7-8750H CPU (2.20GHz / 2.21GHz) and 16 GB RAM memory.
As we can see, in general the diameter strategy has the highest computational complexity, which is comparable to the one of the k-means strategy. The diameter strategy is supposed to be the fastest in computing the cutting direction, but because it produces low quality elements with lost of vertices the refinement time becomes comparable with that of the k-means strategy. The CNN-enhanced strategy outperforms both of them in terms of number of geometrical entities and computational time. It is worth noticing that in the CVT the CNN-enhanced strategy achieves higher or comparable quality with respect to the k-means strategy (Figure \ref{fig:quality metrics} bottom) in much less time. This is because the CNN is classifying CVT elements as \say{cubes}, therefore reducing the computational cost. This means that the CNN is generalizing properly on new samples, because CVT elements were not included in the training set.
\paragraph{Summary.} The diameter strategy has the simplest implementation but it produces low quality elements, raising the computational cost. The k-means strategy produces the most rounded elements and performs better in unstructured grids. It is robust but has a considerable computational cost. The CNN-enhanced strategy well preservers the mesh structure, producing simple elements of uniform size at a low computational cost. However, it is less effective on unstructured grids.

\section{Applications to VEM and PolyDG method}
In this section we test the effectiveness of the proposed refinement strategies within polyhedral
finite element discretizations. We consider the Virtual Element Method (VEM) \cite{beirao2013basic,beirao2014hitchhiker,beirao2016virtual,da2016mixed,da2017high,dassi2018exploring} and the Polygonal Discontinuous Galerkin (PolyDG) method \cite{hesthaven2007nodal,bassi2012flexibility,antonietti2013hp,cangiani2014hp,antonietti2016review,cangiani2017hp} to solve a standard Laplacian problem in 3D: find $u \in H^1_0(\Omega)$ such that
\begin{equation}
 \label{eq:Poisson problem}
\int_\Omega \nabla u \cdot \nabla v = \int_\Omega f v \quad \forall v \in H^1_0(\Omega),
\end{equation}
with $f \in L^2(\Omega)$ a given forcing term. {\color{black} We take $\Omega = (0,1)^3$ and consider a uniform refinement process in Section \ref{section unif} and an adaptive a priori mesh adaptation procedure in Section \ref{section adapt}. In both cases, given a grid of the domain $\Omega$ we compute numerically the solution of problem \eqref{eq:Poisson problem} using either the VEM or the PolyDG method. We then compute the error, in the VEM case in the discrete $H^1_0$-like norm \cite{beirao2014hitchhiker}, while in the PolyDG case in the DG norm 
$$\|v\|_{\text{DG}}^{2}=\sum_{P}\|\nabla v\|_{L^{2}(P)}^{2}+ \sum_F \|\gamma^{1 / 2} \llbracket v \rrbracket \|_{L^{2}\left(F\right)}^{2},$$
where $P$ is a polyhedral mesh element and $F$ is a polygonal element face \cite{arnold2002unified,cockburn2012discontinuous}. Here $\gamma = \alpha p^2/C_{el}$ is the stabilization function, where $\alpha$ is the penalty parameter, $p$ is the polynomial degree used for the approximation, and $C_{el}$ is a function of the shape of the element and is chosen as in \cite{cangiani2014hp}. The jump operator $\llbracket \cdot \rrbracket$, which measures the discontinuity of $v$ across elements, is defined as in \cite{arnold2002unified}.}

\subsection{Uniform refinement} \label{section unif}
In this Section we consider a uniform refinement process, i.e., at each refinement step each mesh element is refined. The forcing term $f$ in \eqref{eq:Poisson problem} is selected in such a way that the exact solution is given by
$$u(x,y) = \sin(\pi x)\sin(\pi y)\sin(\pi z).$$
The grids obtained after three steps of uniform refinement are those already reported in Figure~\ref{fig:corse grids refined} (second to fourth column).\\
We now consider the PolyDG method applied to the grid of cubes for different values of the penalty parameter $\alpha$, using polynomials of degree 1. Larger values of $\alpha$ penalize more discontinuities of the solution, rather than its gradient. In Figure \ref{fig:DG penalty} we show the computed errors as a function of the number of degrees of freedom.
\begin{table}
        \hspace{-0.3\linewidth}
        \begin{tabular}{cM{0.48\linewidth}M{0.48\linewidth}M{0.48\linewidth}}
        
        &  \textbf{ \large $\alpha = 2$} &  \textbf{ \large $\alpha = 5$} &  \textbf{\large $\alpha = 10$} \\
        &
        \includegraphics[width = \linewidth]{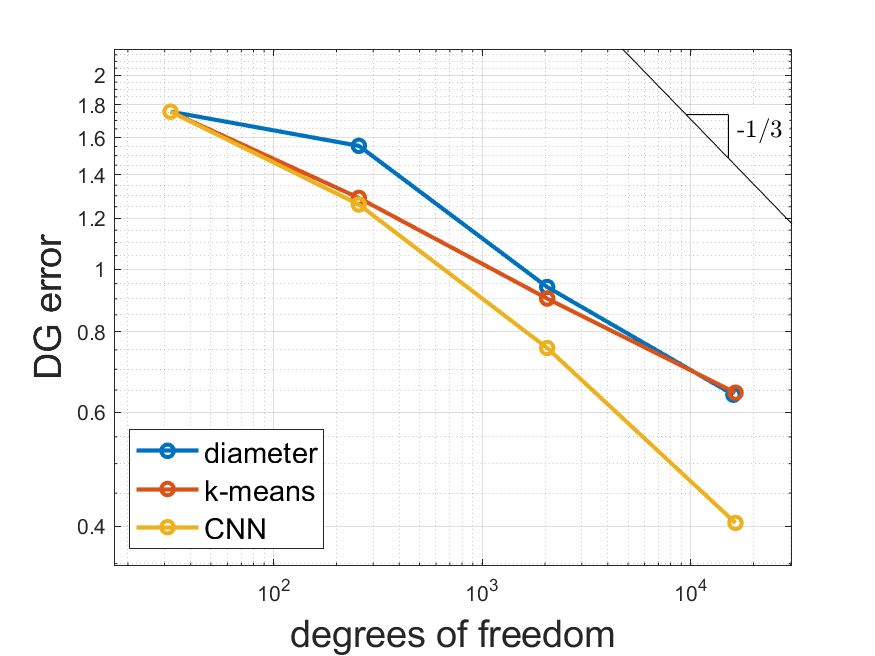} &
        \includegraphics[width = \linewidth]{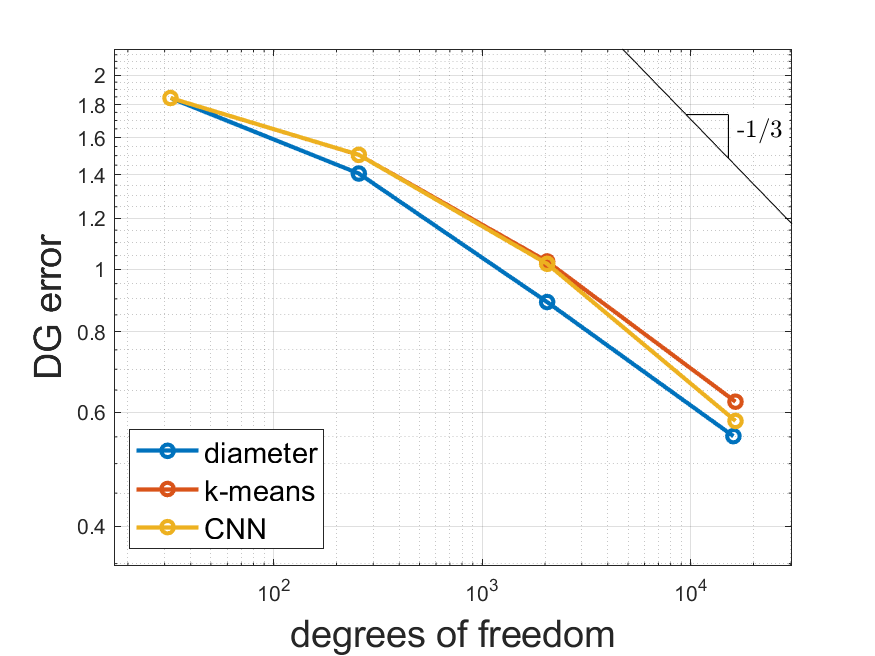} &
        \includegraphics[width = \linewidth]{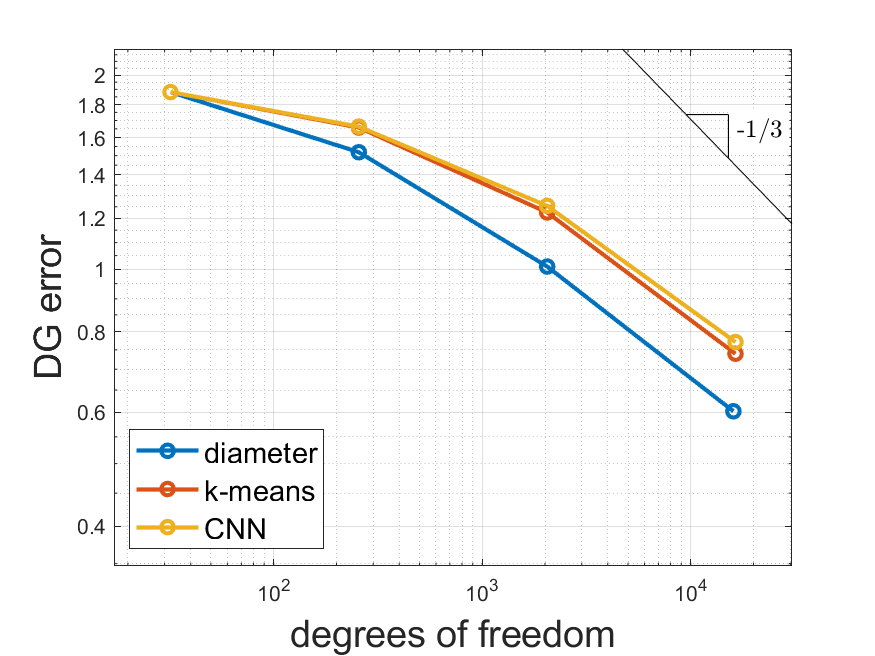}
        \end{tabular}
        \captionsetup{width=1.5\linewidth}
        \captionof{figure}{Uniform refinement test case of Section \ref{section unif}. Computed errors as a function of the number of degrees of freedom. The PolyDG method is applied to the grid of cubes for different values of the penalty parameter $\alpha = 2,5,10$, using polynomials of degree 1. The choice of the most effective refinement strategy (diameter, k-means, CNN) depends on the value of $\alpha$.}
        \label{fig:DG penalty}
\end{table}
As we can see, the CNN-enhance strategy is more effective for $\alpha = 2$, while the diameter strategy for $\alpha = 10$. In order to guarantee convergence, $\alpha$ needs to be larger than a threshold that depends on the shape of the elements. For example, $\alpha = 2$ works for the CNN-enhanced strategy, because it preserves the cube structure of the grid, but not for the other strategies. Since computing the threshold is in general too complex, in the following we will use $\alpha = 5$.\\
In Figure \ref{fig:DG VEM unif} we show the computed errors as a function of the number of degrees of freedom both for the VEM and the PolyDG method, for all grids reported in Figure~\ref{fig:corse grids refined}, using \say{polynomials} of dgree 1.
\begin{table}
    \vspace{-3cm}
    \begin{tabular}{cM{0.425\linewidth}M{0.425\linewidth}}
        
    & \large \textbf{VEM} & \large \textbf{PolyDG} \\
        
    \rotatebox[origin=c]{90}{\large \textbf{tetrahdera} } & 
    \includegraphics[width = \linewidth]{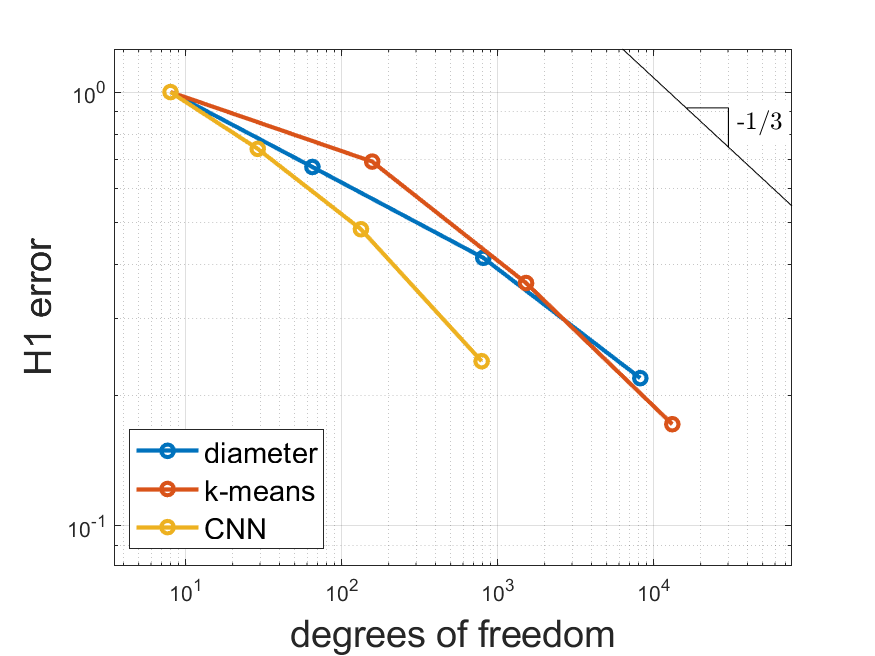} & 
    \includegraphics[width = \linewidth]{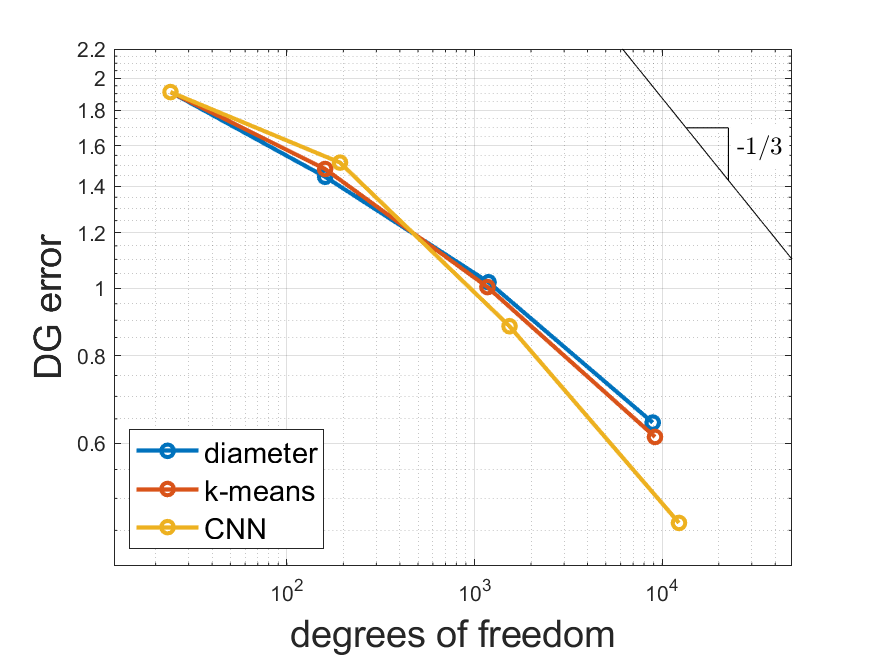}\\
    
    \rotatebox[origin=c]{90}{\large \textbf{cubes} } & 
    \includegraphics[width = \linewidth]{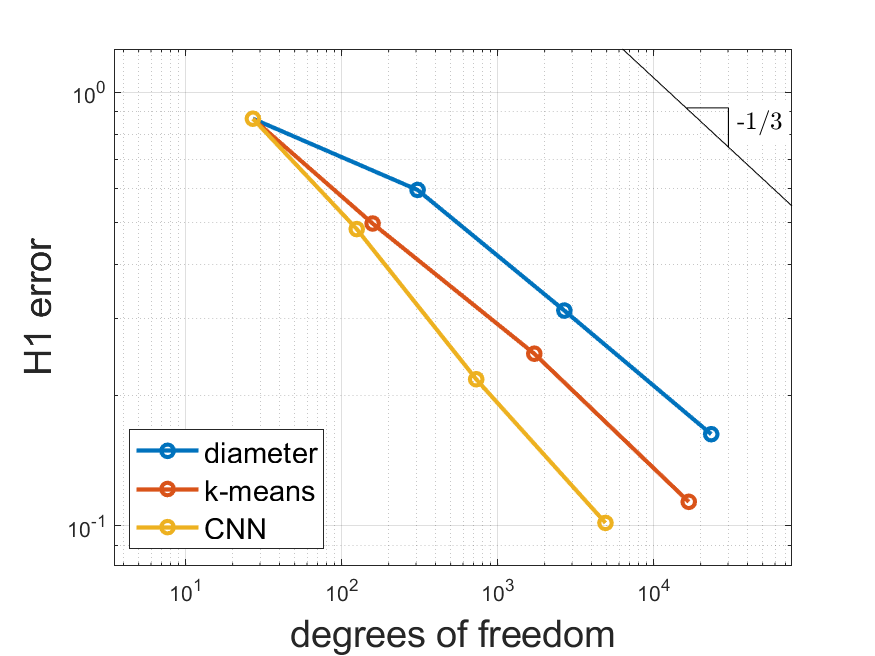} & 
    \includegraphics[width = \linewidth]{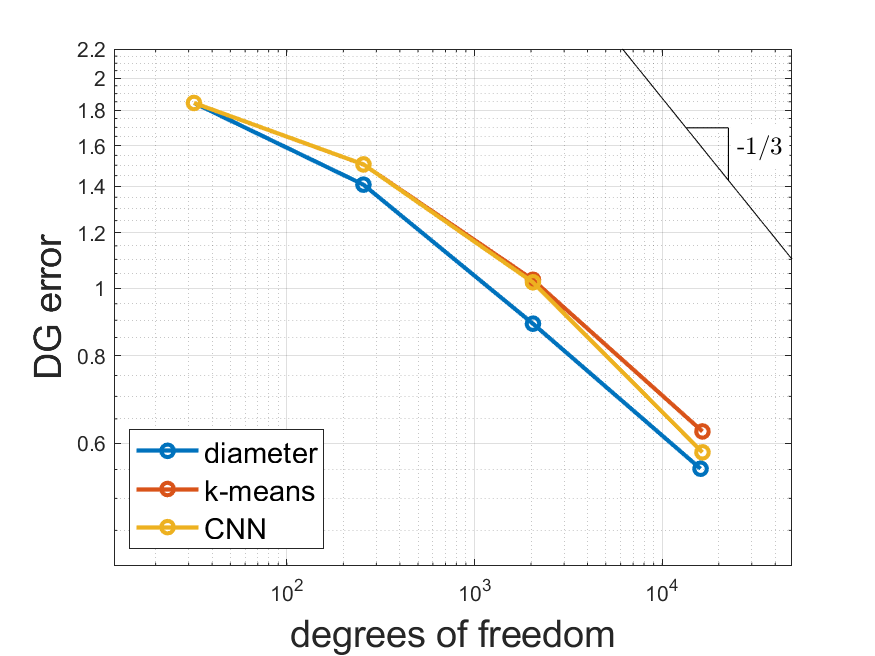}\\
    
    \rotatebox[origin=c]{90}{\large \textbf{prism} } & 
    \includegraphics[width = \linewidth]{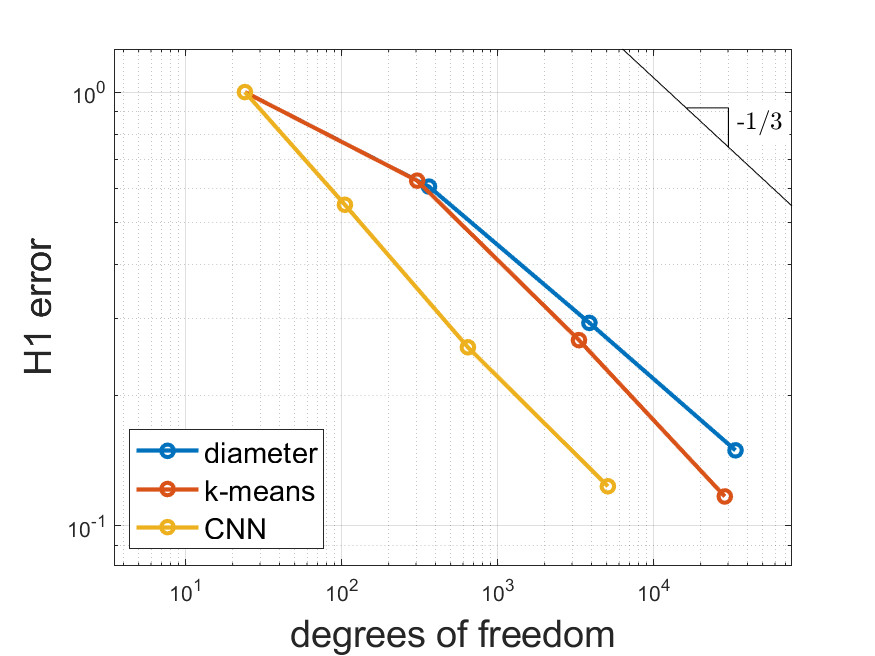} &
    \includegraphics[width = \linewidth]{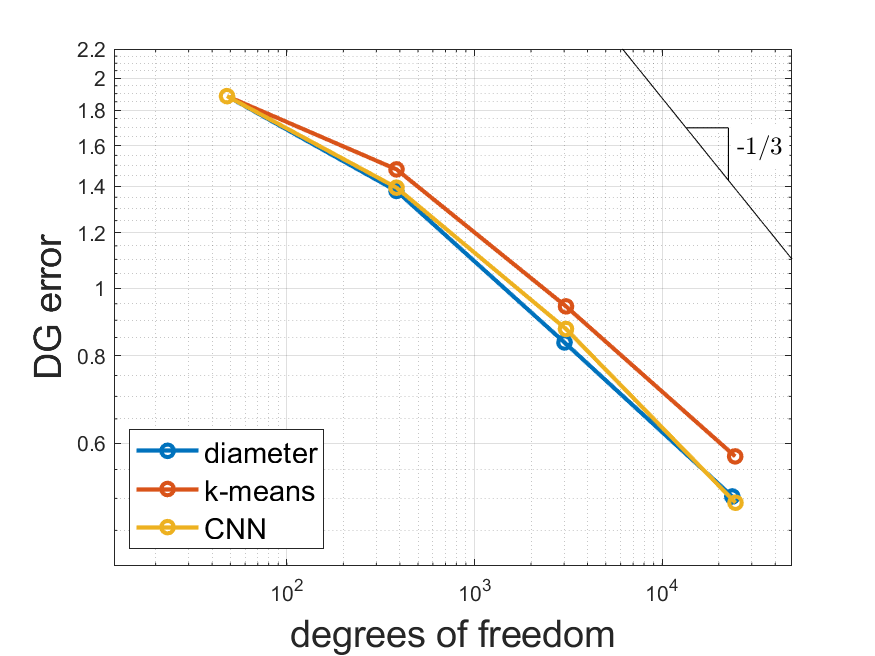}\\
    
    \rotatebox[origin=c]{90}{\large \textbf{Voronoi} } & 
    \includegraphics[width = \linewidth]{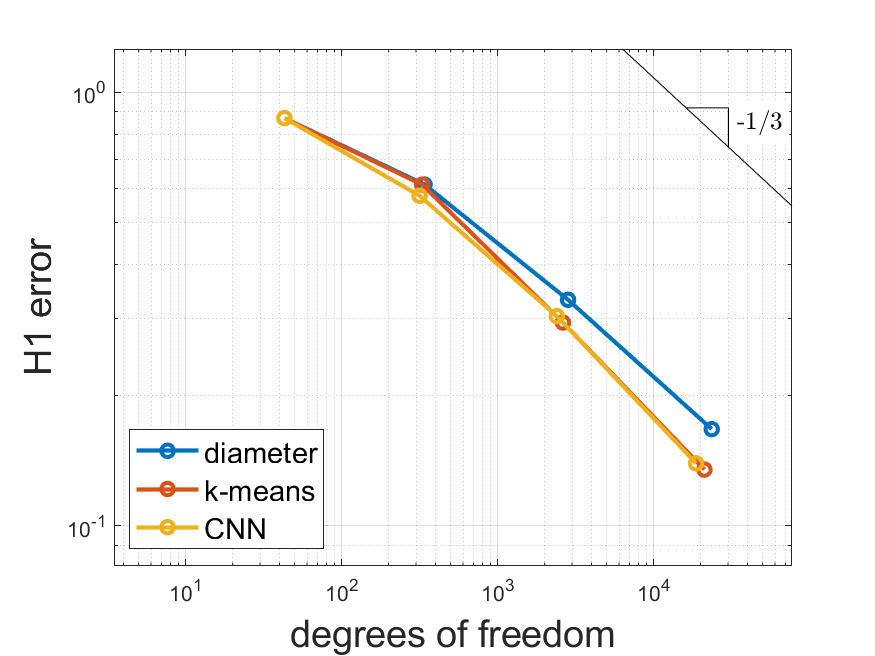}&
    \includegraphics[width = \linewidth]{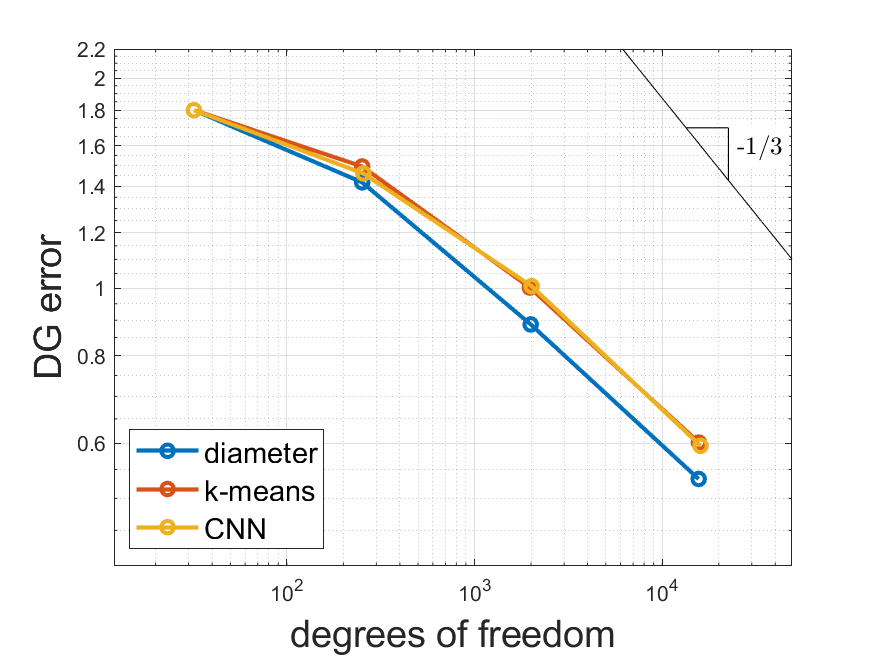}\\
    
    \rotatebox[origin=c]{90}{\large \textbf{CVT} } & 
    \includegraphics[width = \linewidth]{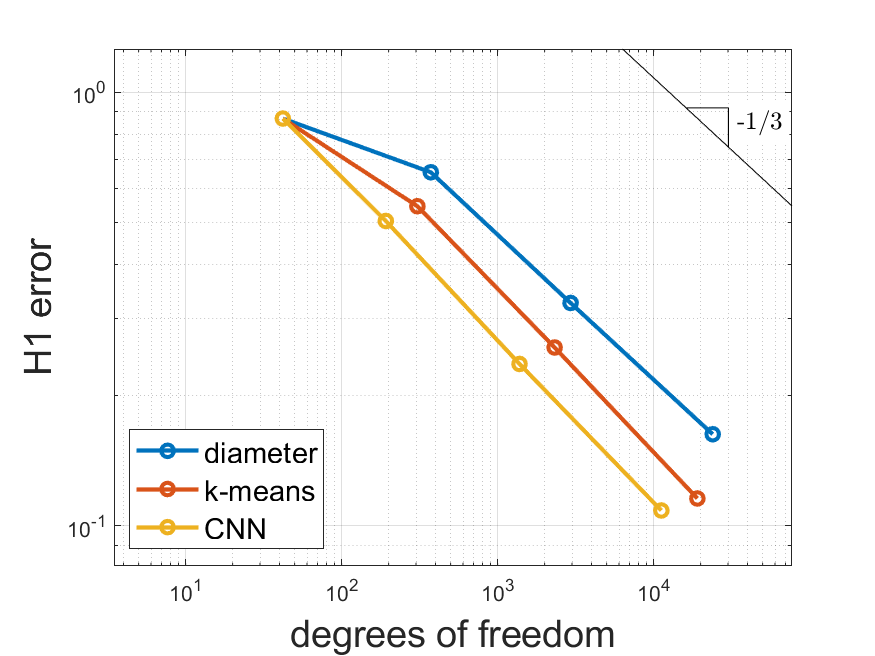}&
    \includegraphics[width = \linewidth]{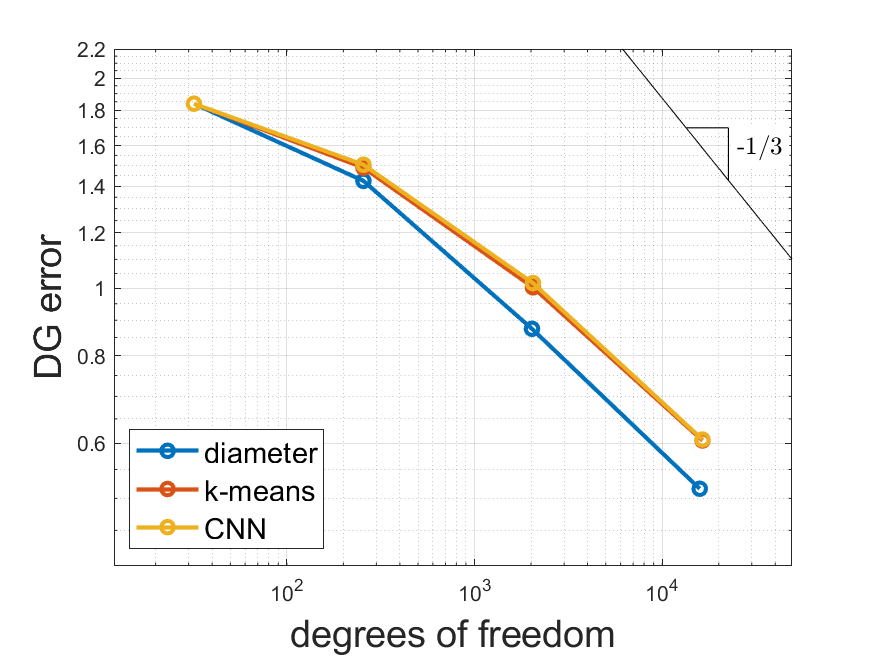}\\
    \end{tabular}
    \captionof{figure}{Uniform refinement test case of Section \ref{section unif}. Computed errors as a function of the number of degrees of freedom. Each row corresponds to the same initial grid refined uniformly with the proposed refinement strategies (diameter, k-means, CNN), while each column corresponds to a different numerical method (VEM left and PolyDG right).}
    \label{fig:DG VEM unif}
\end{table}
When using the PoyDG method, the CNN-enhanced strategy is the most effective in the tetrahedra grid, while the diameter strategy in all other case. To maximize the performance, one may design a CNN-enhanced strategy that classifies polyhedra into two classes only, namely \say{tetrahedron} and \say{other}, in order to apply the ``classical" strategy in former case and the diameter strategy in the latter. 
However, considering that the performance of the different strategies are overall comparable and also that these results may vary depending on the choice of $\alpha$, further analysis would be required in order to effectively exploit the use of ML techniques when employing the PolyDG method. On the contrary, in the VEM case the CNN-enhanced strategy significantly outperforms the others, followed by the k-means strategy. At each step the error obtained with the CNN-enhanced refinement is associated with a lower number of degrees of freedom. Indeed, the error lines are shifted to the left. This sensitivity to mesh distortions may be due to the fact that the VEM is a hybrid method, with unknowns on the elements boundary. In Figure \ref{fig:VEM dgr 2-3 unif} we show that the same results hold in the VEM case also for approximation orders 2 and 3.
\begin{table}
    \vspace{-3cm}
    \hspace{-0.75cm}
    \begin{tabular}{cM{0.425\linewidth}M{0.425\linewidth}}
        
     \large \textbf{VEM}  & \large \textbf{order 2} & \large \textbf{order 3} \\
        
    \rotatebox[origin=c]{90}{\large \textbf{tetrahdera} } & 
    \includegraphics[width = \linewidth]{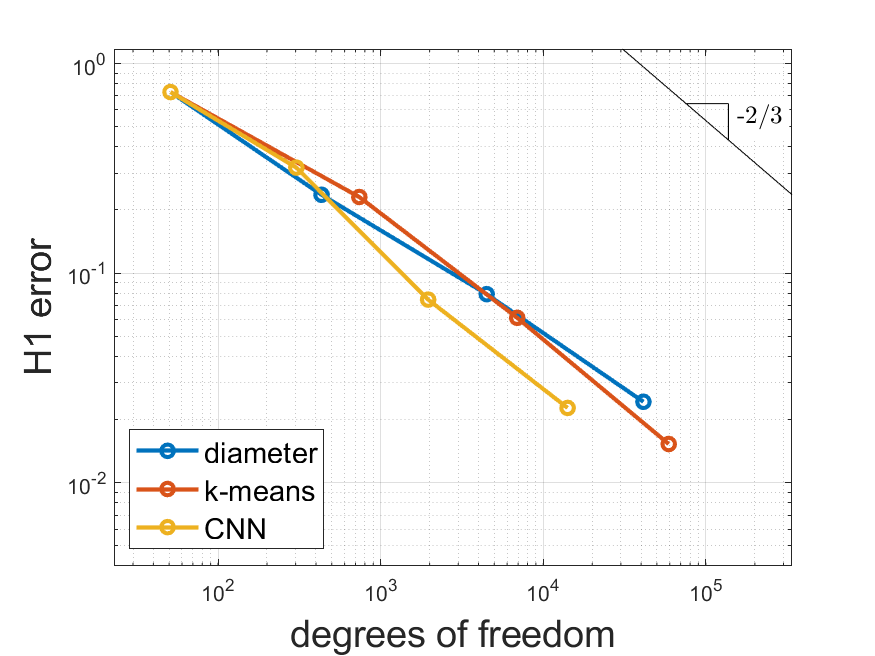} & 
    \includegraphics[width = \linewidth]{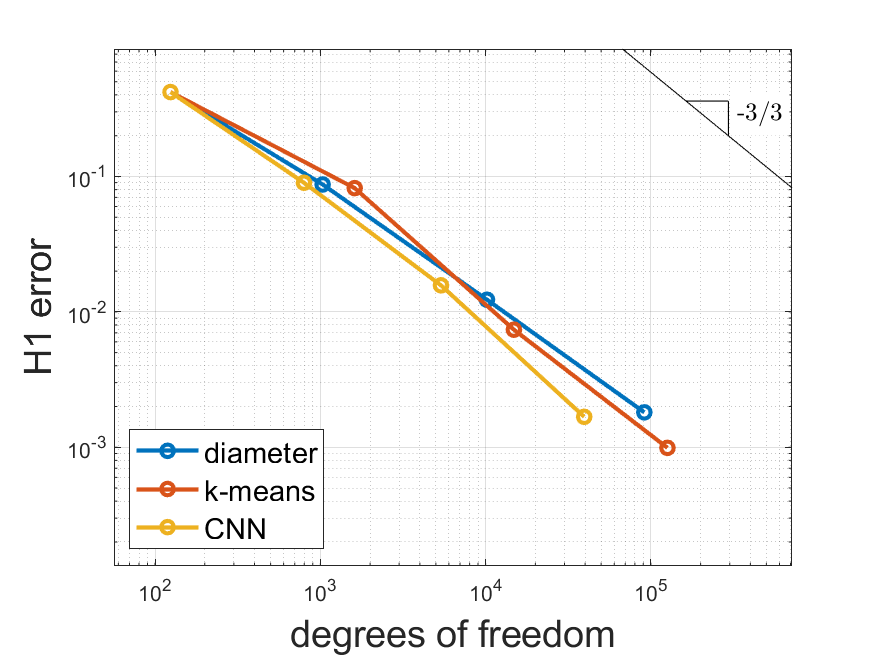}\\
    
    \rotatebox[origin=c]{90}{\large \textbf{cubes} } & 
    \includegraphics[width = \linewidth]{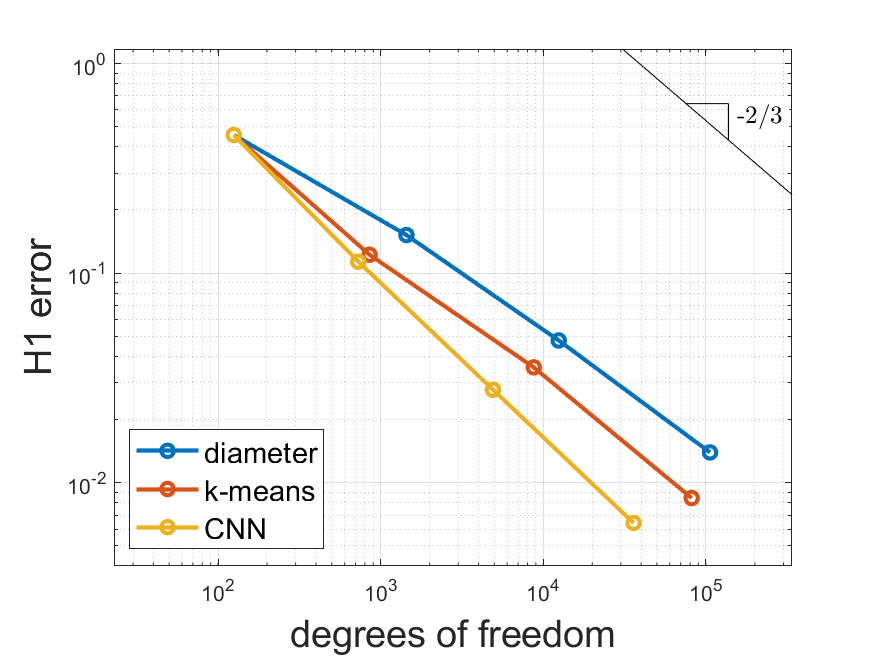} & 
    \includegraphics[width = \linewidth]{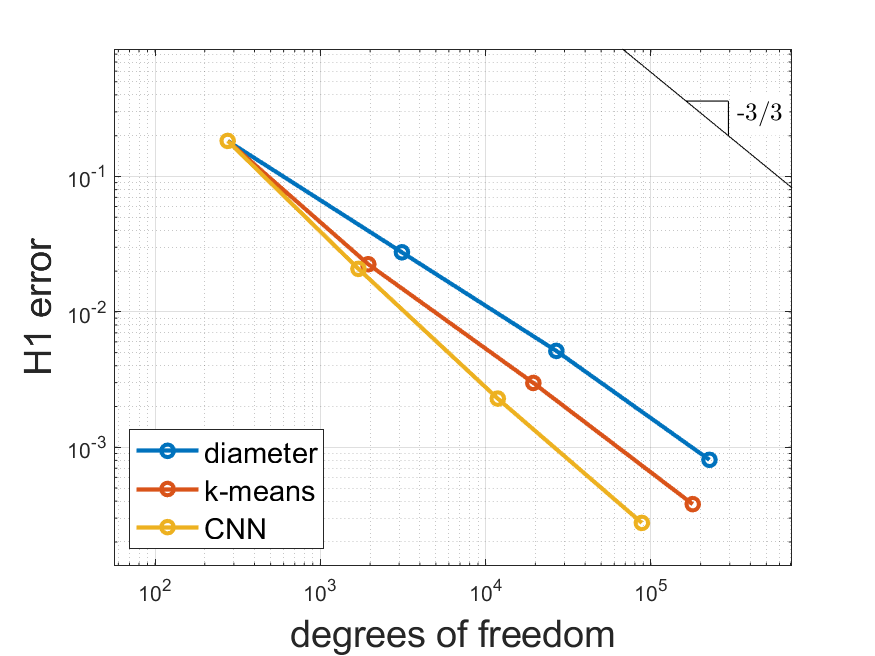}\\
    
    \rotatebox[origin=c]{90}{\large \textbf{prism} } & 
    \includegraphics[width = \linewidth]{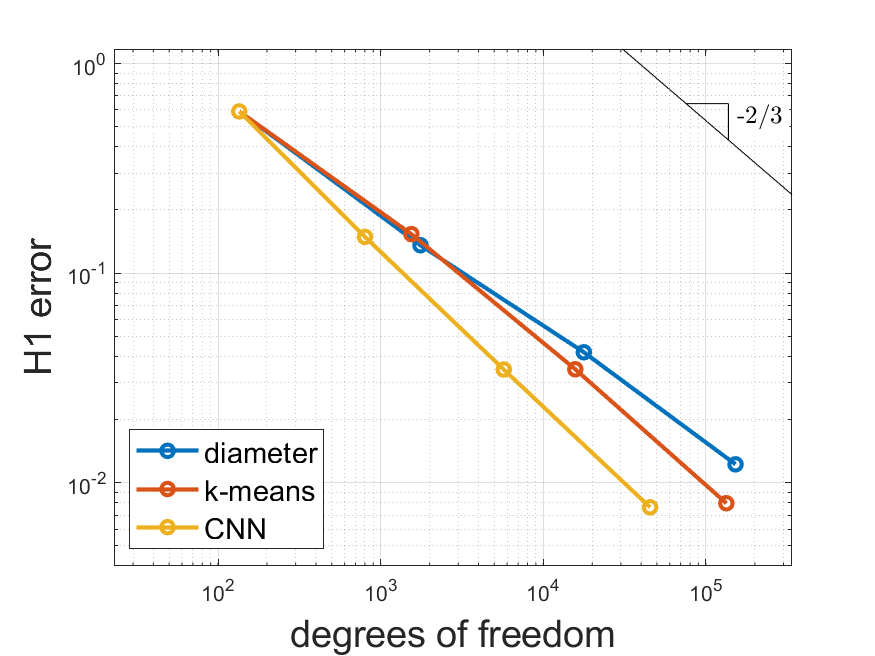} &
    \includegraphics[width = \linewidth]{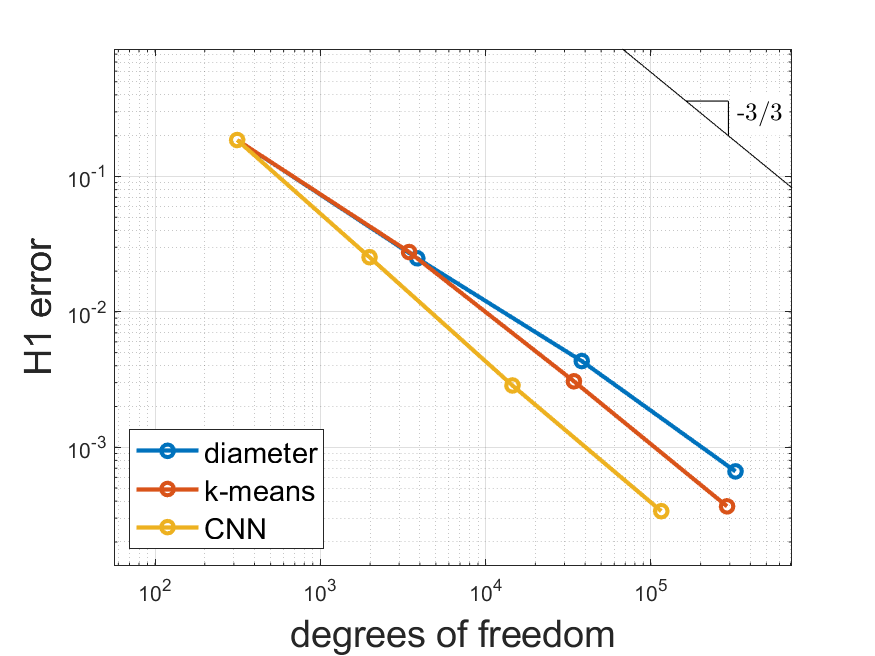}\\
    
    \rotatebox[origin=c]{90}{\large \textbf{Voronoi} } & 
    \includegraphics[width = \linewidth]{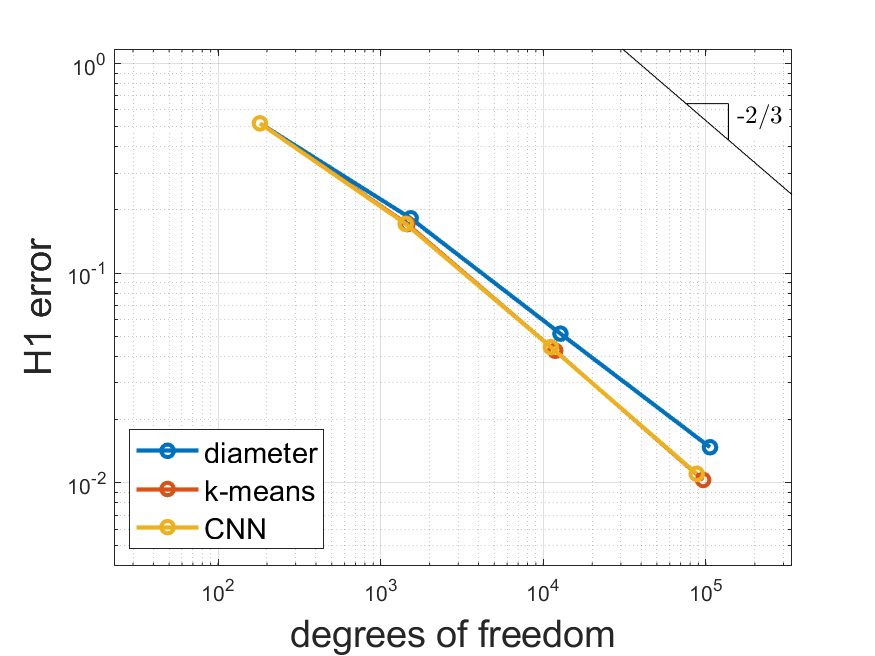}&
    \includegraphics[width = \linewidth]{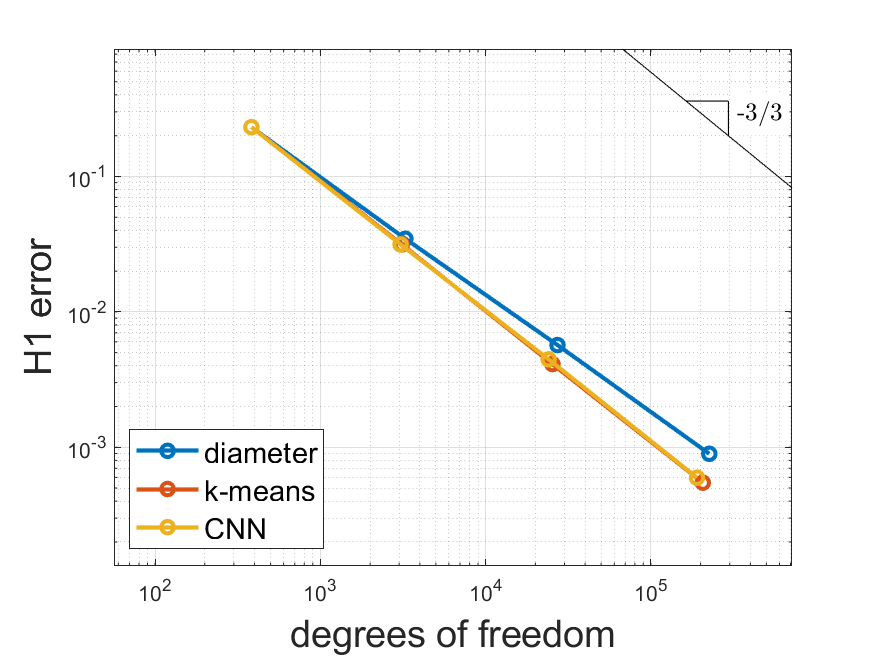}\\
    
    \rotatebox[origin=c]{90}{\large \textbf{CVT} } & 
    \includegraphics[width = \linewidth]{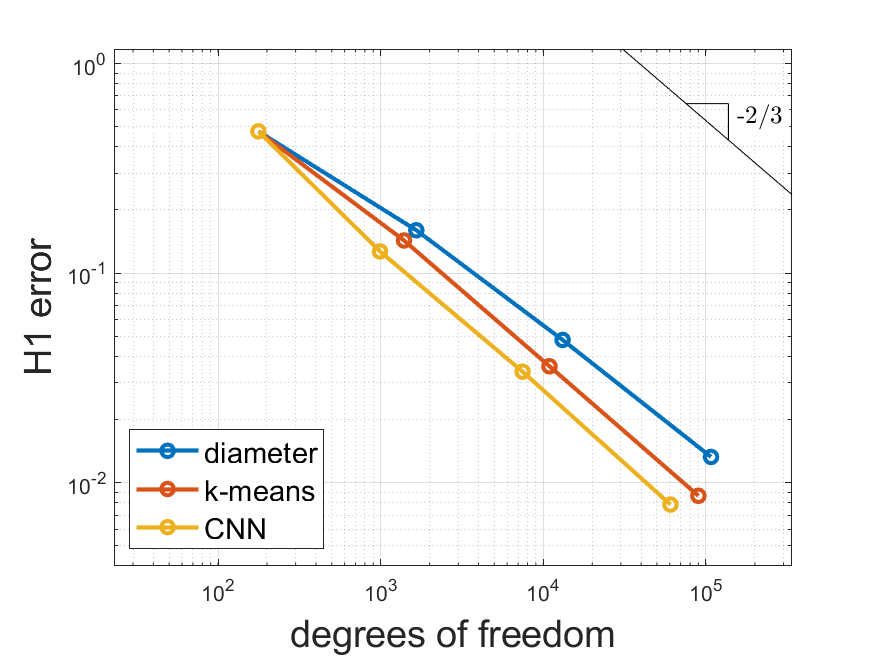}&
    \includegraphics[width = \linewidth]{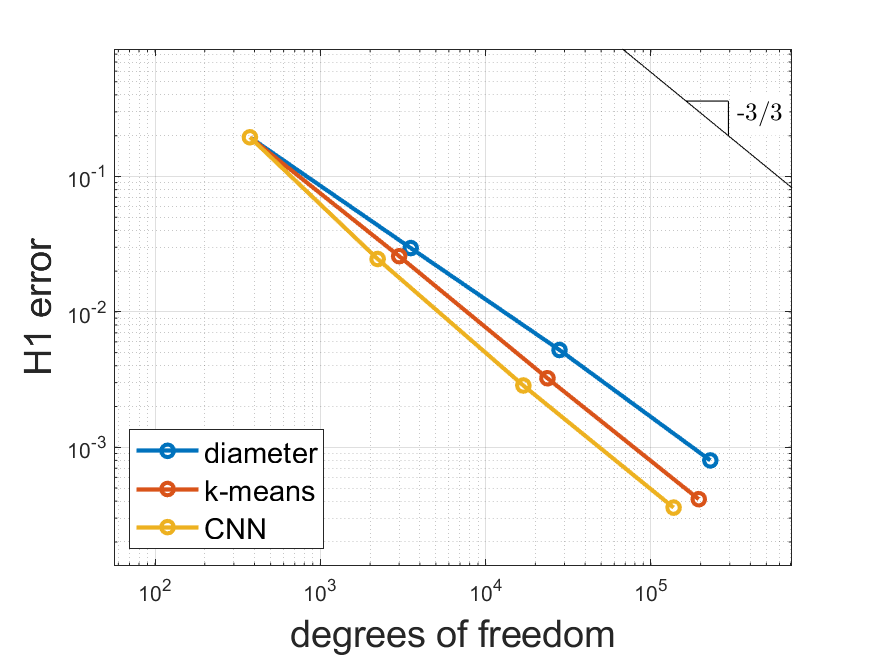}\\
    \end{tabular}
    \captionof{figure}{Uniform refinement test case of Section \ref{section unif}. Computed errors as a function of the number of degrees of freedom. Each row corresponds to the same initial grid refined uniformly with the proposed refinement strategies (diameter, k-means, CNN), while column corresponds to employing the VEM of order 2 and 3.}
    \label{fig:VEM dgr 2-3 unif}
\end{table}

\subsection{Adaptive refinement} \label{section adapt}
\begin{figure}
    \centering
    \includegraphics[width = \linewidth]{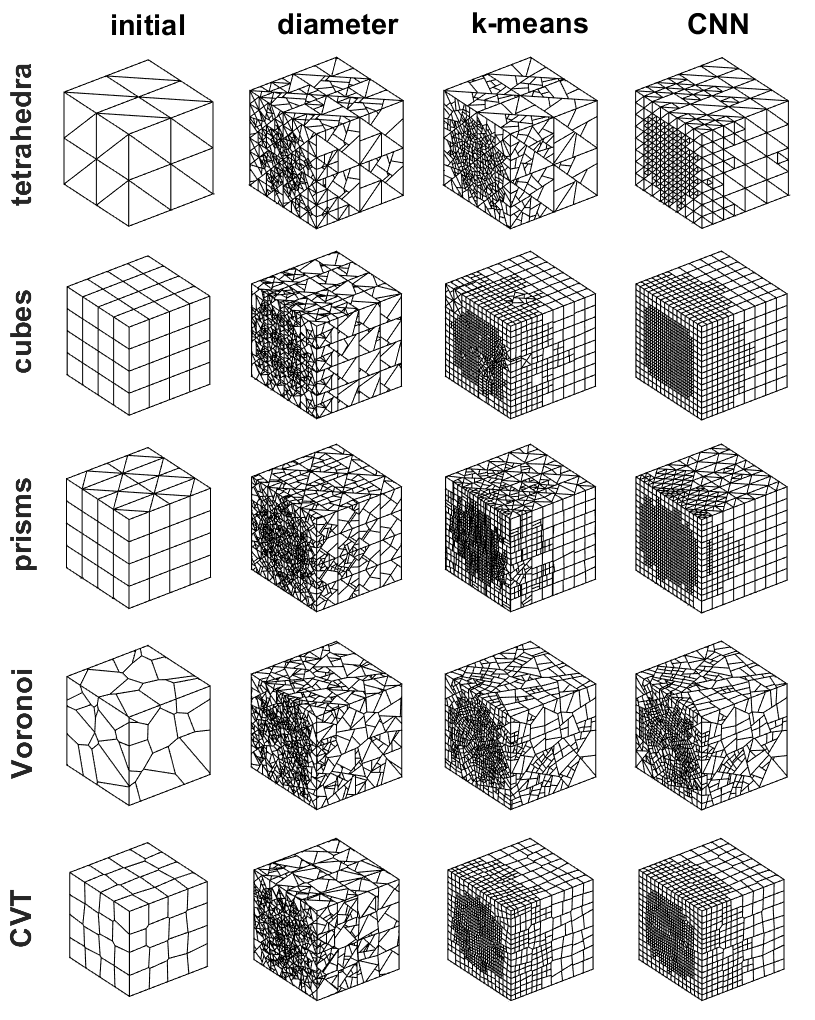}
    \caption{Adaptive refinement test case of Section \ref{section adapt}. Each row corresponds to the same initial grid (tetrahedra, cubes, prisms, Voronoi, CVT), while each column corresponds to the same refinement strategy (diameter, k-means, CNN). Three steps of adaptive refinement have been performed, with a fixed fraction refinement criterion of $r = 0.4$.}
    \label{fig:grids adaptive}
\end{figure}

{\color{black}
In this section we consider an adaptive refinement process. In particular, we proceed as follows:
\begin{enumerate}
    \item Given a grid of the domain $\Omega$ we compute numerically the solution of problem \eqref{eq:Poisson problem} using either the VEM method or the PolyDG.
    \item We compute the error with respect to the exact solution, in the VEM case in the discrete $H^1_0$-like norm \cite{beirao2014hitchhiker}, while in the PolyDG case in the DG norm \cite{cangiani2014hp}. Notice that we did not employ any a posteriori estimator of the error, since we would like to investigate the effect of the proposed refinement strategies.  
    \item We refine a fixed fraction $r$ of elements with the highest error. To refine the marked elements we employ one of the proposed strategies (diameter, k-means, CNN).
\end{enumerate}
}
The forcing term $f$ in \eqref{eq:Poisson problem} is selected in such a way that the exact solution is given by
$$u(x,y) = (1 - e^{-10x})(x-1)\sin(\pi y) \sin(\pi z),$$
that exhibits a boundary layer along  $x = 0$.\\
We take as initial grids the first ones of the previous example and we repeat Steps 1-3 for four times using refinement ratio $r = 0.4$. The adapted meshes employing the PolyDG method are shown in Figure~\ref{fig:grids adaptive}, similar grids are obtained with VEM.
In Figures \ref{fig:DG VEM adapt} and \ref{fig:VEM dgr 2-3 adapt} we show the computed errors for the PolyDG method of degree 1 and the VEM of orders 1, 2 and 3. Results are analogous to the uniform refinement case: in the VEM case there is a clear advantage in using the CNN-enhanced refinement strategy. Also in this case at each step the error obtained with the CNN-enhanced refinement is associated with a lower number of degrees of freedom: the error lines are shifted to the left. In the PolyDG case all strategies have comparable performance. 
\begin{table}
    \vspace{-3cm}
    \begin{tabular}{cM{0.425\linewidth}M{0.425\linewidth}}
        
    & \large \textbf{VEM} & \large \textbf{PolyDG} \\
        
    \rotatebox[origin=c]{90}{\large \textbf{tetrahdera} } & 
    \includegraphics[width = \linewidth]{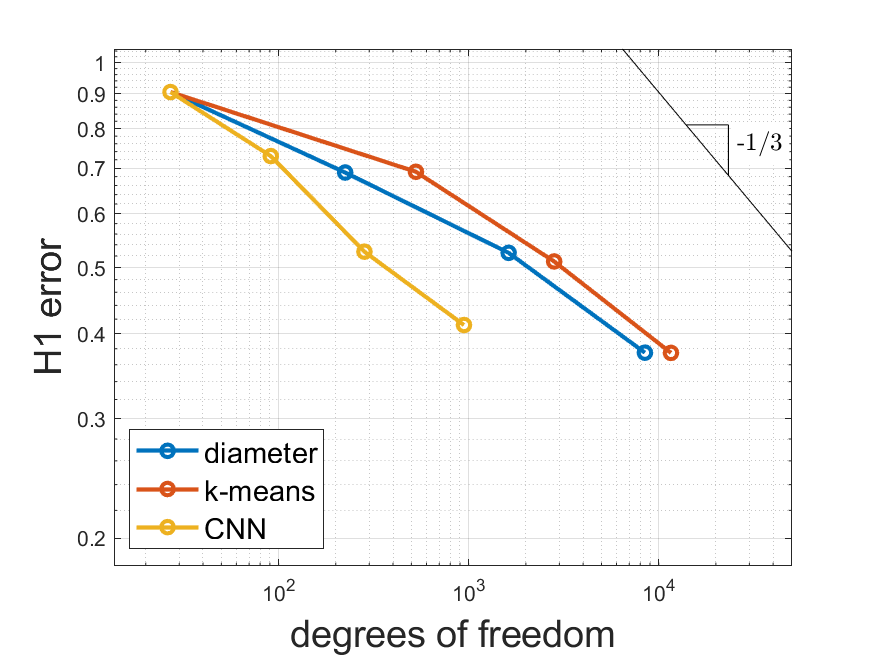} & 
    \includegraphics[width = \linewidth]{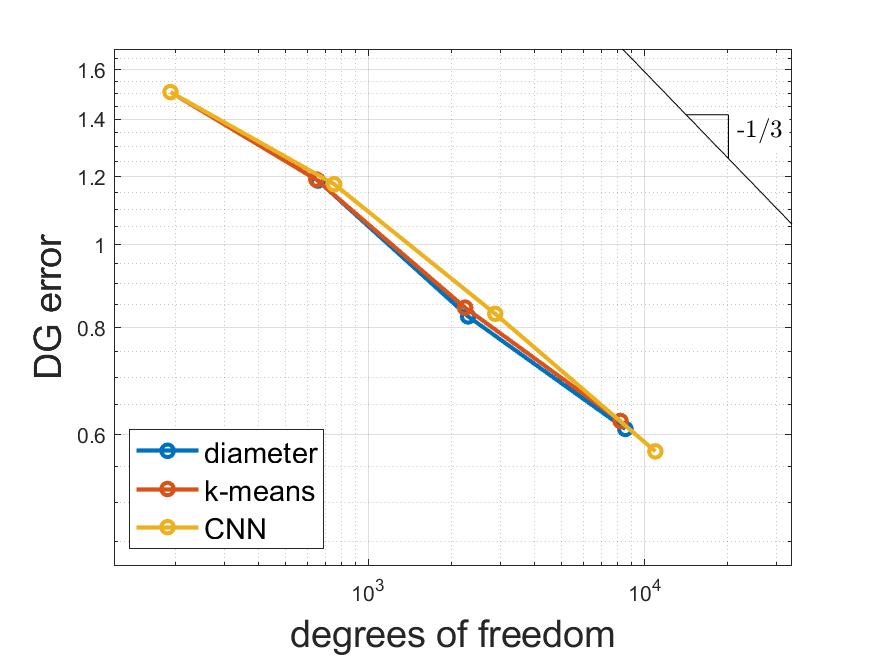}\\
    
    \rotatebox[origin=c]{90}{\large \textbf{cubes} } & 
    \includegraphics[width = \linewidth]{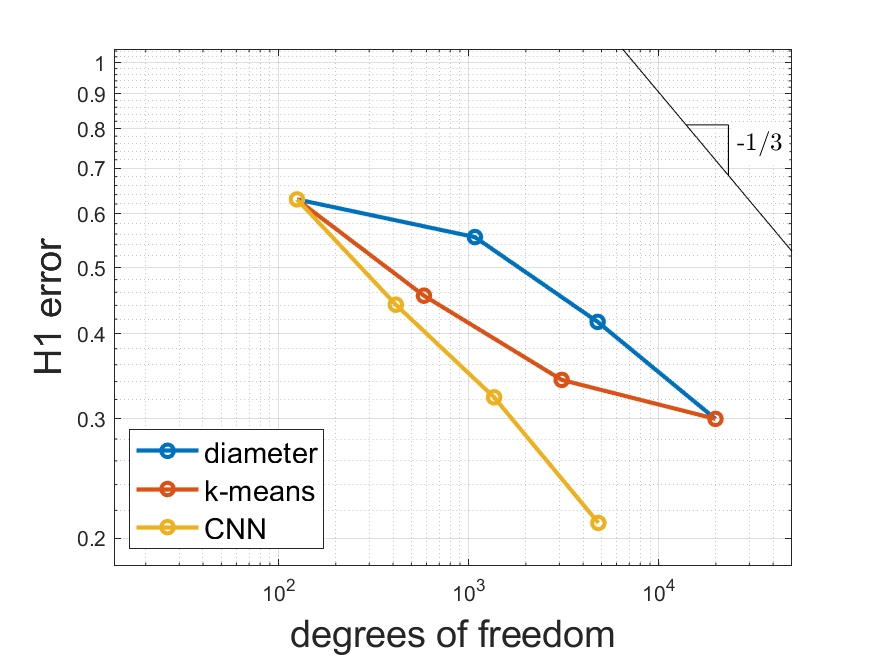} & 
    \includegraphics[width = \linewidth]{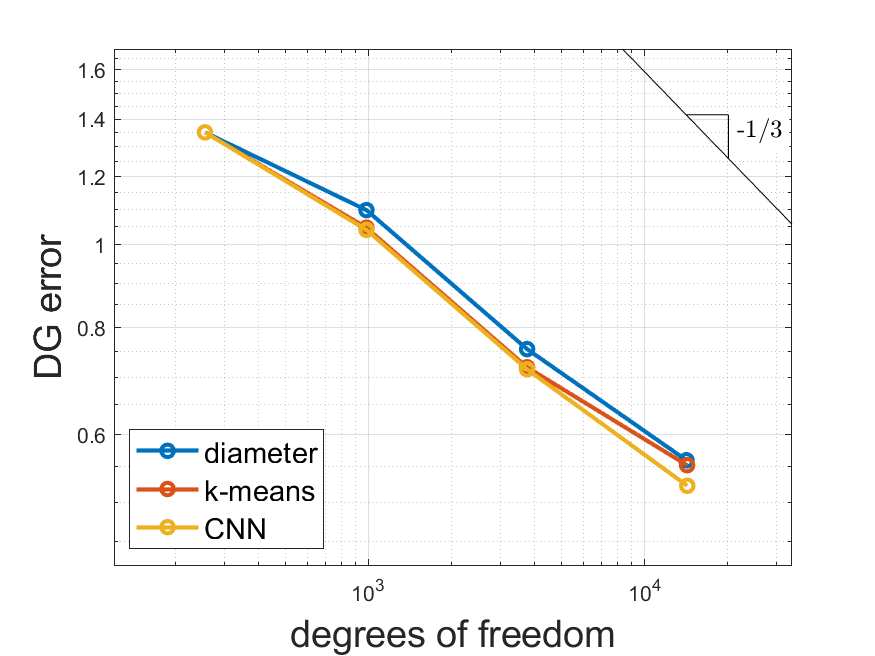}\\
    
    \rotatebox[origin=c]{90}{\large \textbf{prism} } & 
    \includegraphics[width = \linewidth]{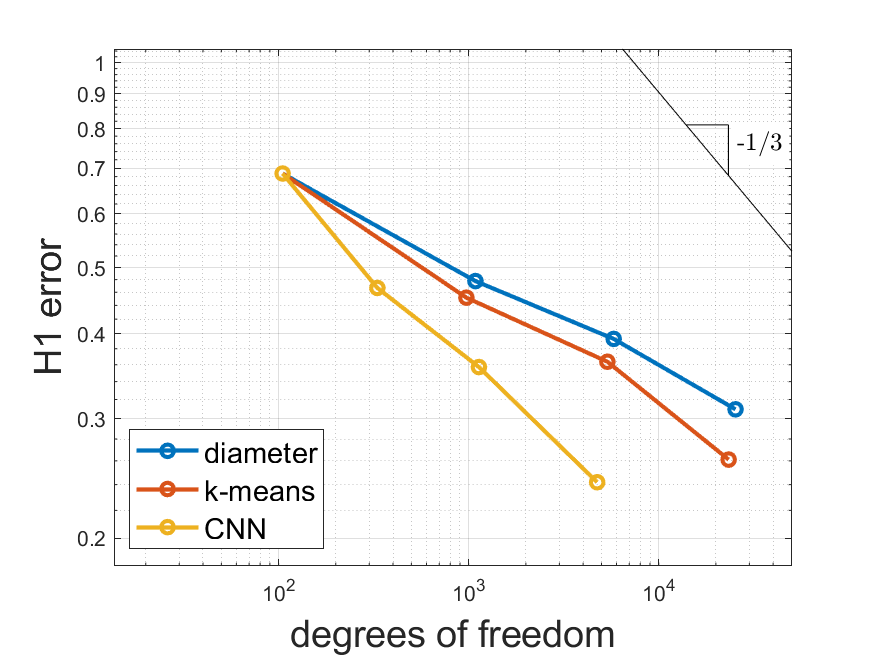} &
    \includegraphics[width = \linewidth]{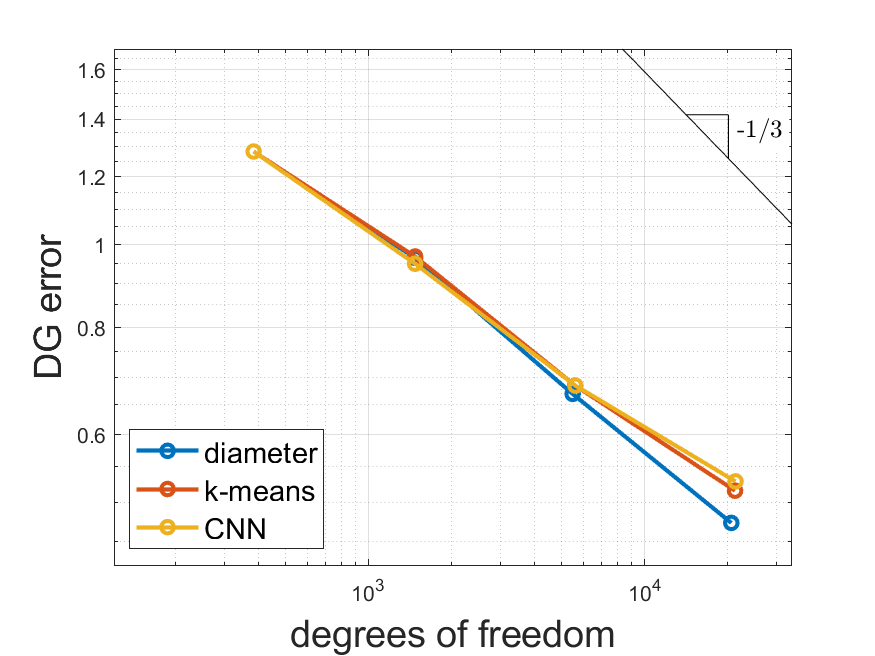}\\
    
    \rotatebox[origin=c]{90}{\large \textbf{Voronoi} } & 
    \includegraphics[width = \linewidth]{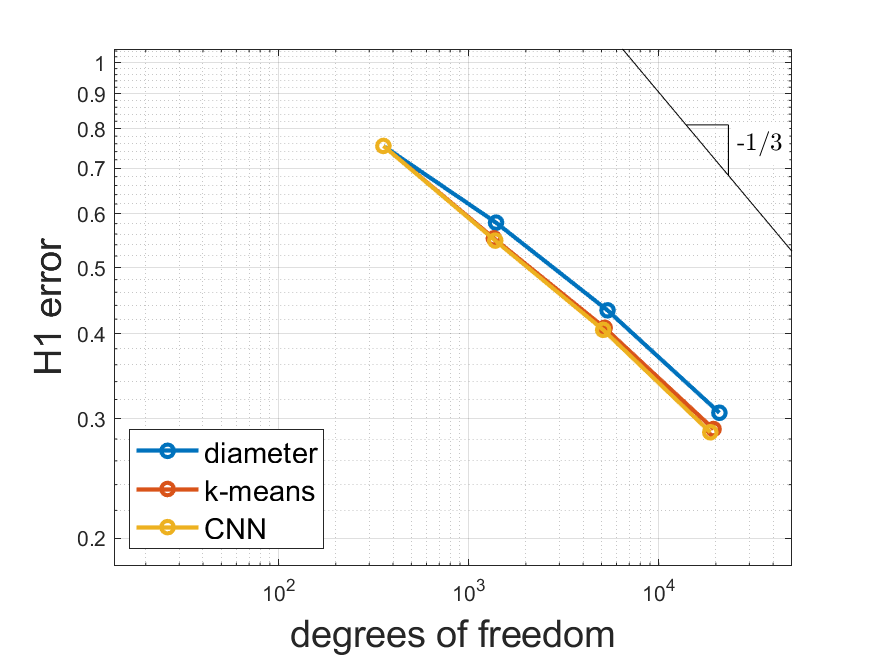}&
    \includegraphics[width = \linewidth]{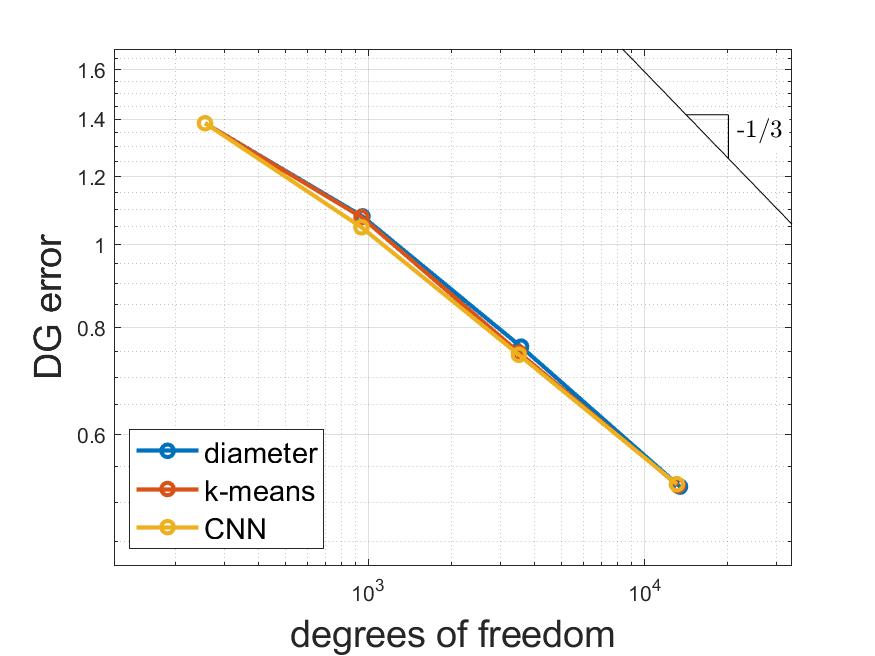}\\
    
    \rotatebox[origin=c]{90}{\large \textbf{CVT} } & 
    \includegraphics[width = \linewidth]{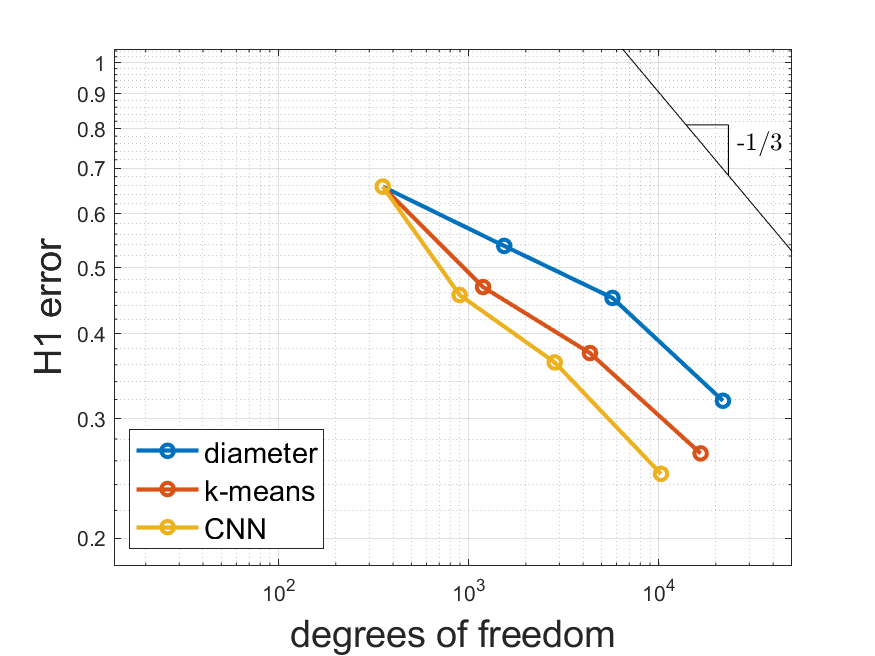}&
    \includegraphics[width = \linewidth]{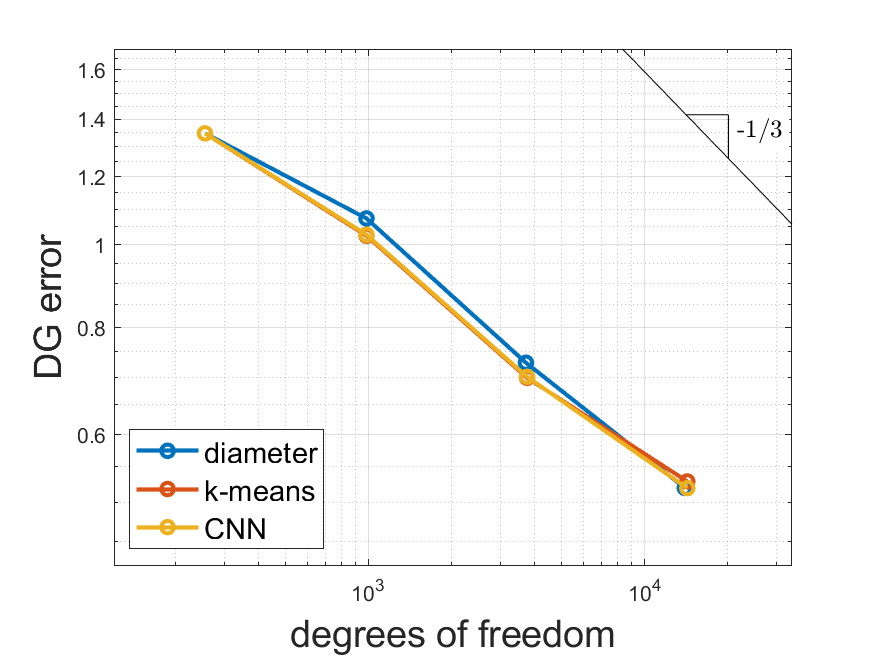}\\
    \end{tabular}
    \captionof{figure}{Adaptive refinement test case of Section \ref{section unif}. Computed errors as a function of the number of degrees of freedom. Each row corresponds to the same initial grid refined uniformly with the proposed refinement strategies (diameter, k-means, CNN), while each column corresponds to a different numerical method (VEM left and PolyDG right).}
    \label{fig:DG VEM adapt}
\end{table}

\begin{table}
    \vspace{-3cm}
    \hspace{-0.75cm}
    \begin{tabular}{cM{0.425\linewidth}M{0.425\linewidth}}
        
     \large \textbf{VEM}  & \large \textbf{order 2} & \large \textbf{order 3} \\
        
    \rotatebox[origin=c]{90}{\large \textbf{tetrahdera} } & 
    \includegraphics[width = \linewidth]{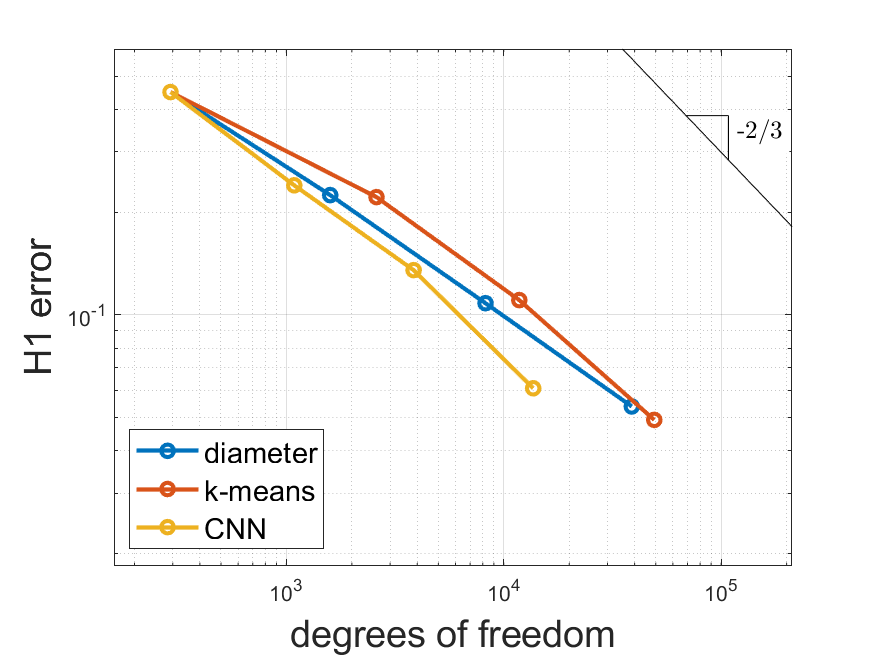} & 
    \includegraphics[width = \linewidth]{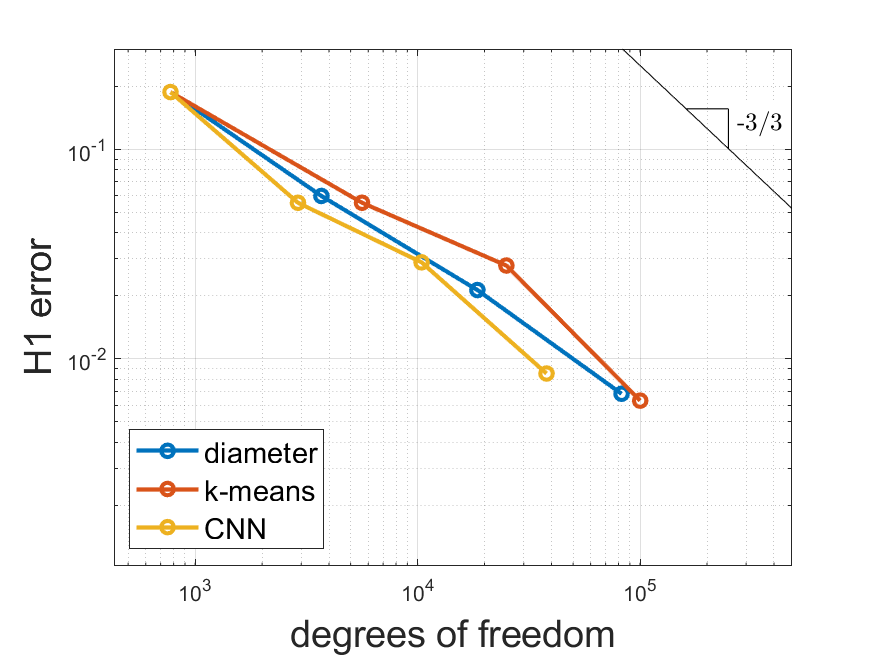}\\
    
    \rotatebox[origin=c]{90}{\large \textbf{cubes} } & 
    \includegraphics[width = \linewidth]{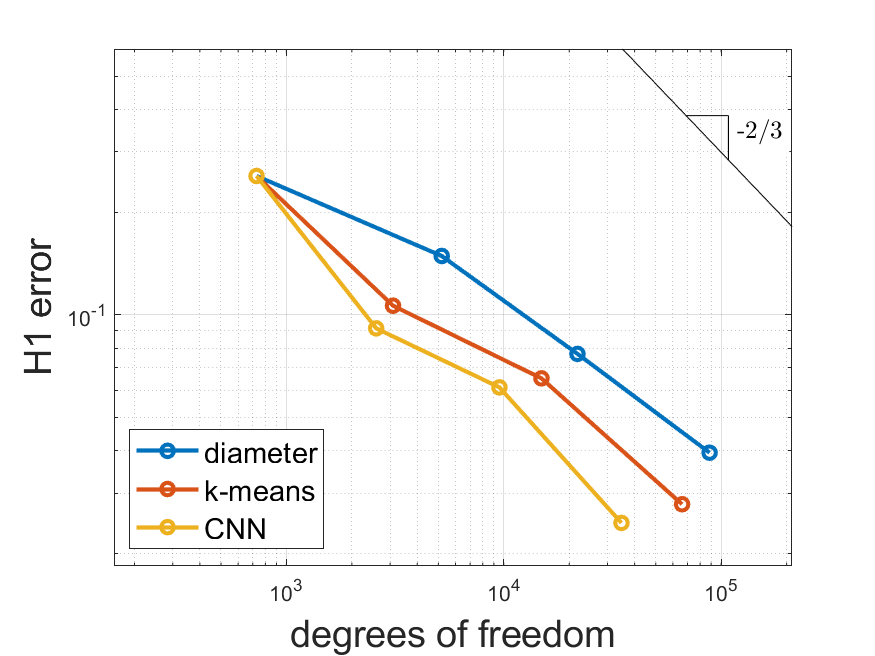} & 
    \includegraphics[width = \linewidth]{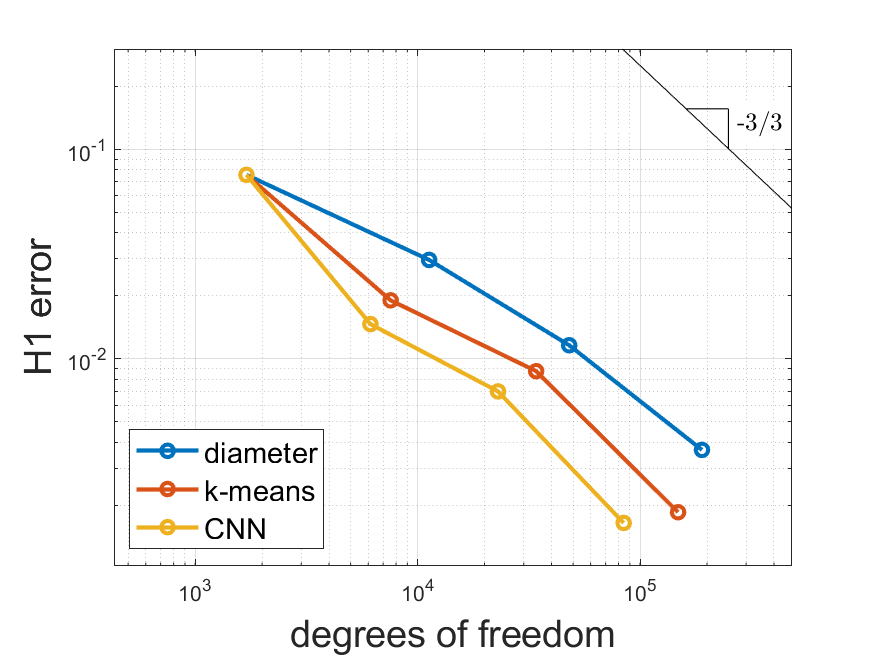}\\
    
    \rotatebox[origin=c]{90}{\large \textbf{prism} } & 
    \includegraphics[width = \linewidth]{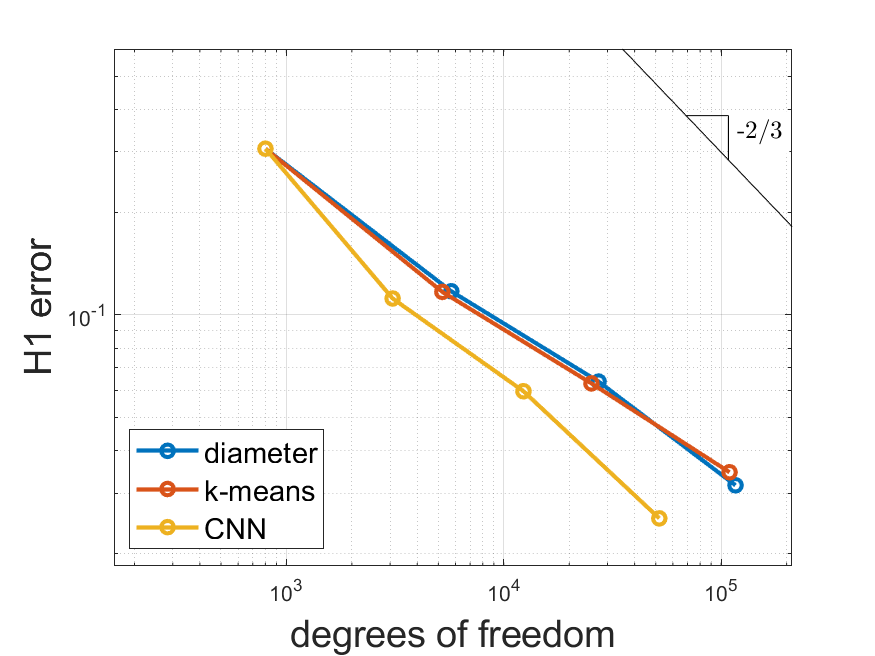} &
    \includegraphics[width = \linewidth]{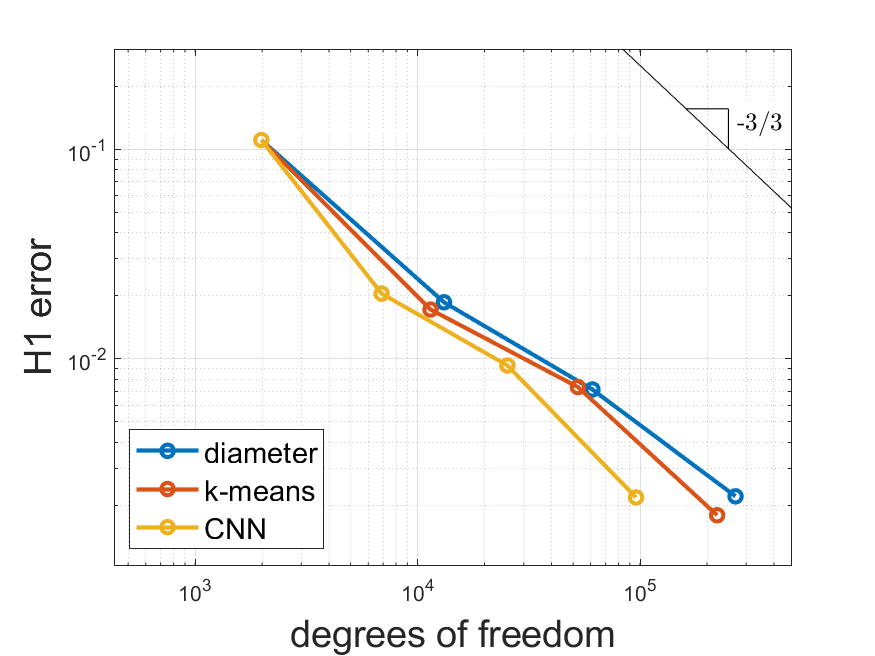}\\
    
    \rotatebox[origin=c]{90}{\large \textbf{Voronoi} } & 
    \includegraphics[width = \linewidth]{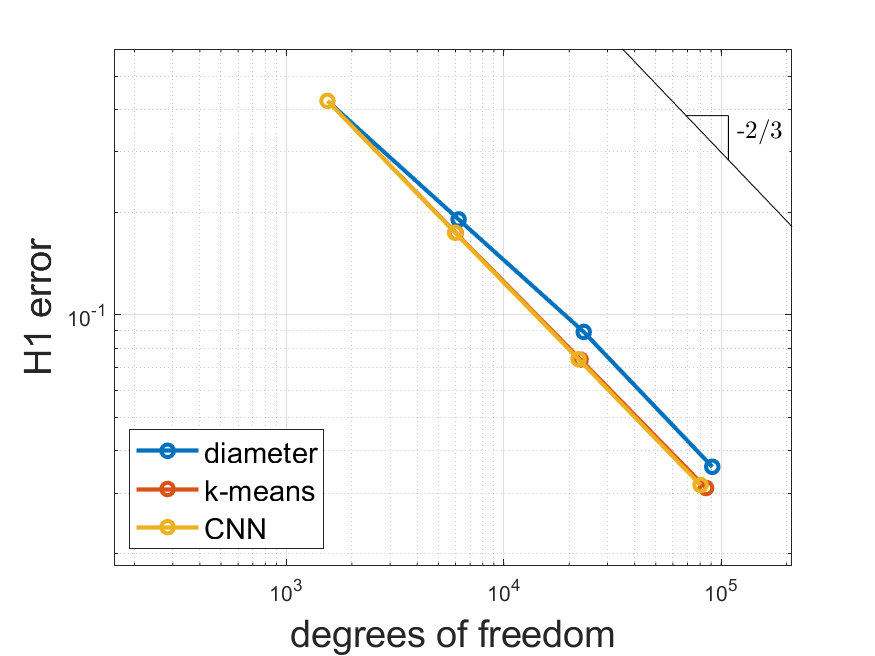}&
    \includegraphics[width = \linewidth]{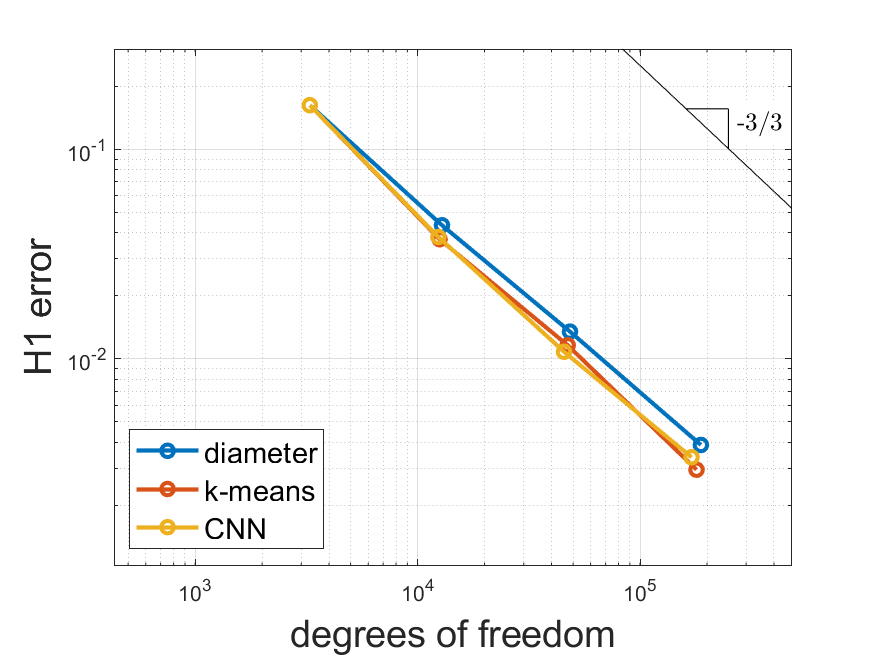}\\
    
    \rotatebox[origin=c]{90}{\large \textbf{CVT} } & 
    \includegraphics[width = \linewidth]{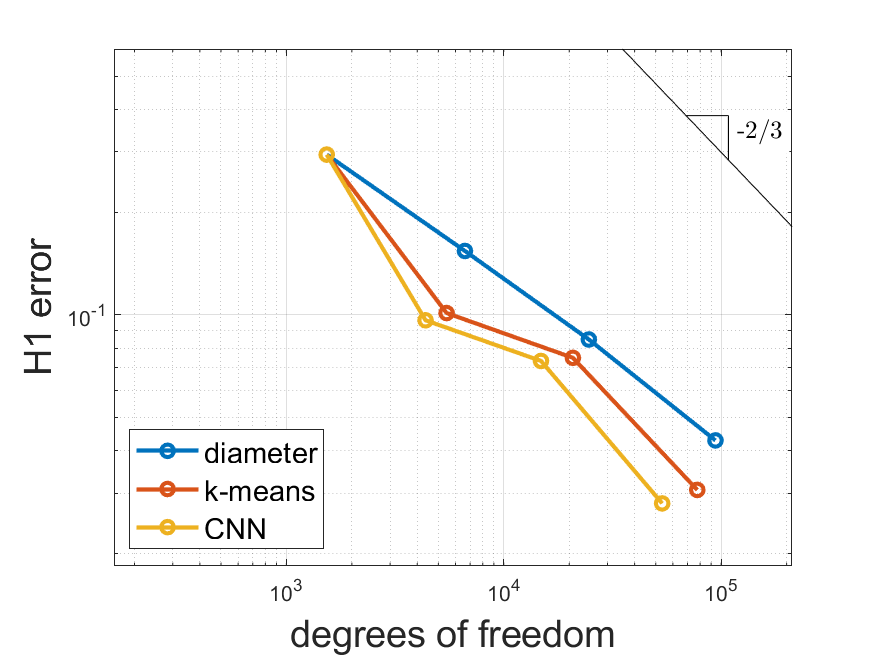}&
    \includegraphics[width = \linewidth]{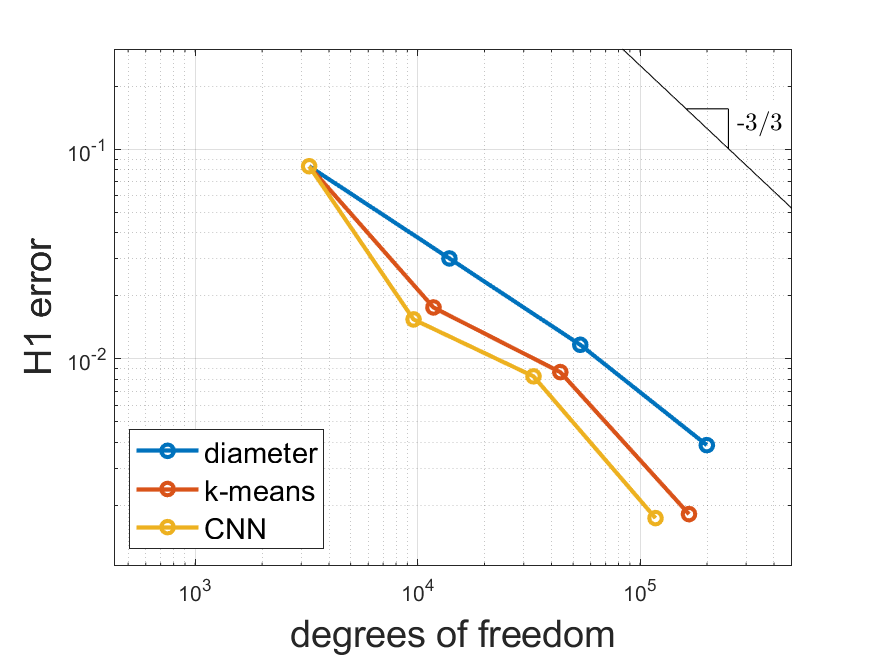}\\
    \end{tabular}
    \captionof{figure}{Adaptive refinement test case of Section \ref{section unif}. Computed errors as a function of the number of degrees of freedom. Each row corresponds to the same initial grid refined uniformly with the proposed refinement strategies (diameter, k-means, CNN), while column corresponds to employing the VEM of order 2 and 3.}
    \label{fig:VEM dgr 2-3 adapt}
\end{table}
\section{Conclusions}
~{\color{black}In this work, we propose the extension to three dimensions of a paradigm based on ML techniques to enhance existing polygonal grid refinement strategies \cite{ANTONIETTI2022110900}, within polyhedral finite element discretizations of partial differential equations. In particular, the k-means algorithm is used to learn a clustered representation of the element that is used to perform the partition. This strategy is a variation of the well known Centroidal Voronoi Tessellation. It produces shape-regular elements and it is robust when applied on unstructured grids.\\
Another approach consists in using a CNN as a classifier for polyhedra, that associates to each element the most suitable refinement criterion, including the k-means. This approach has the advantage of being modular: any refinement strategy can be employed, as well as any numerical method can be used to solve the differential problem. Moreover, the CNN has low computational cost once trained, which makes it appealing especially in three dimensions where processing mesh elements is much more expensive than in two dimensions. The use of CNNs allows to exploit strategies explicitly based on the “shape” of the polyhedron, which would not be possible otherwise unless explicit and costly geometrical checks are performed. As a result, the initial structure and quality of the grid is preserved, significantly lowering the computational cost and the number of vertices, edges and faces required to refine the mesh. These ML-enhanced refinement strategies are particularly beneficial for the VEM, which is sensitive to mesh distortions and/or vertex proliferation, as the error is decreasing faster using the same number of degrees of freedom. On the other hand, as expected the PolyDG method is less sensitive as the discretization space is not associated with any geometrical entity.}\\
In terms of future research lines, we plan to use neural networks to drive agglomeration strategies, either by using CNNs to classify the shape of polyhedra or by using graph neural networks to classify the connectivity graph of mesh elements. These will be of paramount importance when designing multilevel solvers.

\subsection*{CRediT authorship contribution statement}
\textbf{P.F. Antonietti:} Conceptualization, Funding acquisition, Methodology, Project administration, Resources, Supervision, Writing - review and editing.\\
\textbf{F. Dassi:} Software, Validation, Visualization, Writing - review and editing.\\
\textbf{E. Manuzzi:} Conceptualization, Data curation, Investigation, Methodology, Software, Validation, Visualization, Writing – original draft.

\subsection*{Declaration of competing interest}
The authors declare that they have no known competing financial interests or personal relationships that could have appeared to influence the work reported in this paper.

\subsection*{Acknowledgements}
Funding: P.F. Antonietti has been partially supported by the Ministero dell’Università e della Ricerca [PRIN grant numbers 201744KLJL and 20204LN5N5]. P. F. Antonietti, F. Dassi and E. Manuzzi are members of INDAM-GNCS.\\
We thank L. Beir{\~a}o da Veiga for his support.

\newpage

\printbibliography

\end{document}